\documentclass[11pt]{amsart}      %%-*-latex-*-
\usepackage{amsmath,amsthm}
\usepackage{amssymb}  
\usepackage{latexsym}        
\usepackage{euscript}  
\usepackage{graphics}      
\usepackage[frame,ps,dvips,matrix,arrow,curve,rotate]{xy}

\newtheorem{thm}{Theorem}[section]
\newtheorem*{thma}{THEOREM A}

\newtheorem*{propb}{PROPOSITION B}
\newtheorem*{thmc}{THEOREM C}
\newtheorem*{thmd}{THEOREM D}
\newtheorem{defn}[thm]{Definition}
\newtheorem{prop}[thm]{Proposition}
\newtheorem{claim}[thm]{Claim}
\newtheorem{cor}[thm]{Corollary}
\newtheorem{lemma}[thm]{Lemma}

\theoremstyle{definition}  % Bold headings and Roman body text.
\newtheorem{example}[thm]{Example}
\newtheorem{examples}[thm]{Examples}
\newtheorem{facts}[thm]{Facts}
\newtheorem{fact}[thm]{Fact}
\newtheorem{exercise}[thm]{Exercise}
\newtheorem{note}[thm]{Note}
\newtheorem{remark}[thm]{Remark}
\newtheorem{notation}[thm]{Notation}

\newcommand{\dfn}{\textbf} % Make defined words bold.
\newcommand{\cat}[1]{{\EuScript #1}}

\newcommand{\cE}{\cat{E}}
\newcommand{\cF}{\cat{F}}
\newcommand{\cG}{\cat{G}}

\newcommand{\cM}{\cat{M}}

\newcommand{\cO}{\cat{O}}

\newcommand{\cU}{\cat{U}}
\newcommand{\cV}{\cat{V}}
\newcommand{\cW}{\cat{W}}

\newcommand{\bi}{\begin{itemize}}
\newcommand{\ei}{\end{itemize}}
\newcommand{\ba}{\begin{array}}
\newcommand{\ea}{\end{array}}
\newcommand{\bc}{\begin{center}}
\newcommand{\ec}{\end{center}}
\newcommand{\bt}{\begin{tabular}}
\newcommand{\et}{\end{tabular}}
\newcommand{\beqn}{\begin{equation}}
\newcommand{\eeqn}{\end{equation}}

\newcommand{\bp}{\begin{proof}}
\newcommand{\ep}{\end{proof}}
\newcommand{\blemma}{\begin{lemma}}
\newcommand{\elemma}{\end{lemma}}
\newcommand{\bprop}{\begin{prop}}
\newcommand{\eprop}{\end{prop}}
\newcommand{\bthm}{\begin{thm}}
\newcommand{\ethm}{\end{thm}}
\newcommand{\bremark}{\begin{remark}}
\newcommand{\eremark}{\end{remark}}
\newcommand{\bcor}{\begin{cor}}
\newcommand{\ecor}{\end{cor}}
\newcommand{\bnote}{\begin{note}}
\newcommand{\enote}{\end{note}}
\newcommand{\bdefn}{\begin{defn}}
\newcommand{\edefn}{\end{defn}}
\newcommand{\bexample}{\begin{example}}
\newcommand{\eexample}{\end{example}}
\newcommand{\bclaim}{\begin{claim}}
\newcommand{\eclaim}{\end{claim}}
\newcommand{\bexercise}{\begin{exercise}}
\newcommand{\eexercise}{\end{exercise}}
\newcommand{\bnotation}{\begin{notation}}
\newcommand{\enotation}{\end{notation}}
\newcommand{\bexamples}{\begin{examples}}
\newcommand{\eexamples}{\end{examples}}
\newcommand{\bfacts}{\begin{facts}}
\newcommand{\efacts}{\end{facts}}
\newcommand{\bfact}{\begin{fact}}
\newcommand{\efact}{\end{fact}}

\newcommand{\E}{\textbf{\textit{E}}}

\newcommand{\ho}{HO}

\newcommand{\hm}{HM}

\newcommand{\one}{\boldsymbol{1}}  %{1\hspace{-1.3mm} 1}

\newcommand{\field}[1]  {\mathbb{#1}} % Use blackboard bold for these sets

\newcommand{\N}         {\field{N}}
\newcommand{\Z}         {\field{Z}}
\newcommand{\C}         {\field{C}}

\newcommand{\eps}{\epsilon}
\renewcommand{\aa}{\alpha}
\newcommand{\bb}{\beta}
\renewcommand{\gg}{\gamma}

\newcommand{\dd}{\delta}

\renewcommand{\ll}{\lambda}
\newcommand{\Ll}{\Lambda}
\newcommand{\oo}{\omega}

\renewcommand{\ss}{\sigma}

\renewcommand{\tt}{\theta}

\DeclareMathOperator{\im}{Im}  % This is for `Image'; Im is lready taken.
\DeclareMathOperator{\Spec}{Spec}

\newcommand{\del}{\partial}
\newcommand{\ts}{\otimes}                       % tensor product
\newcommand{\ra}{\rightarrow}                   % right arrow
\newcommand{\lra}{\longrightarrow}              % long right arrow
\newcommand{\la}{\leftarrow}                    % left arrow
\newcommand{\llra}[1]{\stackrel{#1}{\lra}}      % labeled long right arrow
\newcommand{\we}{\llra{\sim}}                   % weak equivalence
\newcommand{\inc}{\hookrightarrow}              % inclusion

\newcommand{\vs}{\vspace{2mm}}

\newcommand{\blank}{-}                          % A hyphen, as in (-)xV
\newcommand{\period}    {{\makebox[0pt][l]{\hspace{2pt} .}}}
\newcommand{\comma}     {{\makebox[0pt][l]{\hspace{2pt} ,}}}

\begin{document}

\title{Equivariant Elliptic Cohomology and Rigidity}
  
\author{Ioanid Rosu}

\begin{abstract}
Equivariant elliptic cohomology with complex coefficients was
defined axiomatically by Ginzburg, Kapranov and Vasserot~\cite{GKV}
and constructed by Grojnowski~\cite{Gr}.  We give an invariant
definition of complex $S^1$-equivariant elliptic cohomology, and use
it to give an entirely cohomological proof of the rigidity theorem
of Witten for the elliptic genus.  We also state and prove a rigidity
theorem for families of elliptic genera.
\end{abstract}

\maketitle

\tableofcontents

\section{Introduction} \label{section-introduction}

The classical level 2 elliptic genus is defined
(see Landweber~\cite{La}, p.56) as the Hirzebruch genus
with exponential series the Jacobi sine\footnote{For a
definition of the Jacobi sine $s(x)$ see the beginning
of Section~\ref{section-pushforwards}.}.
It is intimately related with the mysterious field of 
elliptic cohomology (see Segal~\cite{Seg2}), and
with string theory (see Witten~\cite{Wi1} and~\cite{Wi2}).
A striking property of the elliptic genus is its 
rigidity with respect to group actions.  
This  was conjectured by Ochanine in~\cite{O1}, and by
Witten in~\cite{Wi1}, where he used string theory
arguments to support it.

Rigorous mathematical proofs for the rigidity of the elliptic
genus were soon given by Taubes~\cite{Ta},
Bott \& Taubes~\cite{BT}, and Liu~\cite{Liu}.  While Bott and
Taubes's proof involved the localization
formula in equivariant K-theory, Liu's proof focused on
the modularity properties of the elliptic genus.  The question
remained however whether one could find a direct connection
between the rigidity theorem and elliptic cohomology.

Earlier on, Atiyah and Hirzebruch~\cite{AH} had used 
pushforwards in equivariant K-theory to prove 
the rigidity of the $\hat{A}$-genus for spin manifolds. 
Following this idea, H. Miller~\cite{Mi} interpreted
the equivariant elliptic genus as a pushforward in 
the completed Borel equivariant cohomology,
and posed the problem of developing and using a 
\emph{noncompleted} $S^1$-equivariant elliptic cohomology,
to prove the rigidity theorem.

In 1994 Grojnowski~\cite{Gr} proposed a noncompleted
equivariant elliptic cohomology theory with complex coefficients.
For $G$ a compact connected Lie group he defined $\E^*_G(\blank)$
as a coherent holomorphic sheaf over a certain variety $X_G$
constructed from a given elliptic curve.  Grojnowski also
constructed pushforwards in this theory.
At about the same time and independently, Ginzburg, Kapranov 
and Vasserot~\cite{GKV} gave an axiomatic description of 
equivariant elliptic cohomology.

Given Grojnowski's construction, it seemed natural to try to
use $S^1$-equivariant elliptic cohomology to prove the rigidity
theorem.  In doing so, we noticed that our proof relies
on a generalization of Bott and Taubes' ``transfer 
formula'' (see~\cite{BT}).  This generalization 
turns out to be essentially equivalent to the existence
of a Thom class (or orientation) in $S^1$-equivariant elliptic
cohomology.

We can generalize the results of this paper in several 
directions.  One is to extend the rigidity theorem to families
of elliptic genera, which we do in Theorem~\ref{rigidity-families}.
Another would be to generalize from $G=S^1$ to an arbitrary connected 
compact Lie group, or to replace complex coefficients 
with rational coefficients for all cohomology theories involved.  
Such generalizations will be treated elsewhere.

\vs

\section{Statement of results} \label{section-results}

All the cohomology theories involved in this paper have complex
coefficients.  If $X$ is a finite $S^1$-$CW$ complex,
$H^*_{S^1}(X)$ denotes its Borel $S^1$-equivariant cohomology
with complex coefficients (see Atiyah and Bott~\cite{AB}).
If $X$ is a point $*$, $H^*_{S^1}(*)\cong\C[u]$.

Let $\cE$ be an elliptic curve over $\C$.  Let $X$ be a finite
$S^1$-$CW$ complex, e.g.\ a compact $S^1$-manifold\footnote{A
compact $S^1$-manifold always has an $S^1$-$CW$ complex
structure: see Alday and Puppe~\cite{AP}.}.  Then, following
Grojnowski~\cite{Gr}, we define $\E^*_{S^1}(X)$, the
$S^1$-equivariant elliptic cohomology of $X$.  This is a coherent
analytic sheaf of $\Z_2$-graded algebras over $\cE$.  We alter
his definition slightly, in order to show that the definition
of $\E^*_{S^1}(X)$ depends only on $X$ and the elliptic curve $\cE$.
Let $\aa$ be a point of $\cE$.  We associate a subgroup $H(\aa)$
of $S^1$ as follows: if $\aa$ is a torsion point of $\cE$ of 
exact order $n$, $H(\aa)= \Z_n$; otherwise, $H(\aa)=S^1$.  We define 
$X^\aa=X^{H(\aa)}$, the subspace of $X$ fixed by $H(\aa)$.
Then we will define a sheaf $\E^*_{S^1}(X)$ over $\cE$ whose
stalk at $\aa$ is
  $$\E^*_{S^1}(X)_\aa = H^*_{S^1}(X^\aa) \ts_{\C[u]}\cO_{\C,0} \period$$
Here $\cO_{\C,0}$ represents the local ring of germs of holomorphic
functions at zero on $\C=\Spec\C[u]$.  In particular, the stalk of
$\E^*_{S^1}(X)$ at zero is $H^*_{S^1}(X) \ts_{\C[u]}\cO_{\C,0}$.

\begin{thma} \label{A}
\textit{$\E^*_{S^1}(X)$ only depends on $X$ and the elliptic curve
$\cE$.  It extends to an $S^1$-equivariant cohomology theory with
values in the category of coherent analytic sheaves of $\Z_2$-graded
algebras over $\cE$.}
\end{thma}

If $f: X \ra Y$ is a complex oriented map between compact
$S^1$-manifolds, Grojnowski also defines equivariant elliptic
pushforwards.  They are maps of sheaves of $\cO_\cE$-modules
$f^E_!:\E^*_{S^1}(X)^{[f]} \ra \E^*_{S^1}(Y)$
satisfying properties similar to those of the usual pushforward
(see Dyer~\cite{Dy}).
$\E^*_{S^1}(X)^{[f]}$ has the same stalks as $\E^*_{S^1}(X)$,
but the gluing maps are different.

If $Y$ is a point, then $f^E_!(1)$ on the stalks at zero is the
$S^1$-equivariant elliptic genus of $X$ (which is a power series
in $u$).  By analyzing in detail the construction
of $f^E_!$, we obtain the following interesting result, which
answers a question posed by H. Miller and answered independently
by Dessai and Jung~\cite{DJ}.

\begin{propb} \label{B}
\textit{The $S^1$-equivariant elliptic genus of a compact
$S^1$-manifold is the Taylor expansion at zero of a function on
$\C$ which is holomorphic at zero and meromorphic everywhere.}
\end{propb}

Grojnowski's construction raises a few natural questions.  First,
can we say more about $\E^*_{S^1}(X)^{[f]}$?  The answer is
given in Proposition~\ref{identification-thom-sheaf}, where we
show that, up to an invertible sheaf, $\E^*_{S^1}(X)^{[f]}$
is the $S^1$-equivariant elliptic cohomology of the Thom space
of the stable normal bundle to $f$.  (In fact, if we enlarge
our category of equivariant $CW$-complexes to include 
equivariant spectra, we can show that 
$\E^*_{S^1}(X)^{[f]}$ is the reduced $\E^*_{S^1}$ of a Thom
spectrum $X^{-Tf}$.  See the discussion after 
Proposition~\ref{identification-thom-sheaf} for details.)

This suggests looking for a Thom section (orientation) in 
$\E^*_{S^1}(X)^{[f]}$.  More generally, given a real oriented
vector bundle $V \ra X$, we can twist $\E^*_{S^1}(X)$ in a similar
way to obtain a sheaf, which we denote by $\E^*_{S^1}(X)^{[V]}$.
For the rest of this section we regard all the sheaves not on $\cE$,
but on a double cover $\tilde{\cE}$ of $\cE$.  The reason for
this is given in the beginning of Subsection~\ref{section-rigidity}.2.
So when does a Thom section exist in $\E^*_{S^1}(X)^{[V]}$?
The answer is the following key result.

\begin{thmc} \label{C}
\textit{If $\; V\ra X$ is a spin $S^1$-vector bundle over a finite
$S^1$-$CW$ complex, then the element 1 in the stalk of
$\E^*_{S^1}(X)^{[V]}$ at zero extends to a global section,
called the Thom section.}
\end{thmc}

The proof of Theorem C is essentially a generalization
of Bott and Taubes' ``transfer formula'' (see ~\cite{BT}).
Indeed, when we try to extend 1 to a global section, we
see that the only points where we encounter difficulties
are certain torsion points of $\cE$ which we call \emph{special}
(as defined in the beginnning of Section~\ref{section-ell-defn}).
But extending our section at a special point $\aa$ amounts to
lifting a class from $H^*_{S^1}(X^{S^1}) \ts_{\C[u]}\cO_{\C,0}$
to $H^*_{S^1}(X^\aa) \ts_{\C[u]}\cO_{\C,0}$ via the restriction
map $i^*: H^*_{S^1}(X^\aa) \ts_{\C[u]}\cO_{\C,0} \ra
H^*_{S^1}(X^{S^1}) \ts_{\C[u]}\cO_{\C,0}$.  This is not a
problem, except when we have two different connected components
of $X^{S^1}$ inside one connected component of $X^\aa$.  Then
the two natural lifts differ up to a sign, which can be shown
to disappear if $V$ is spin.  This observation is due to
Bott and Taubes, and is the centerpiece of their ``transfer
formula.''

Given Theorem C, the rigidity theorem of Witten follows easily:
Let $X$ be a compact spin $S^1$-manifold.  Then the 
$S^1$-equivariant pushforward of $f:X\ra *$ is a map of
sheaves $f^E_!:\E^*_{S^1}(X)^{[f]} \ra \E^*_{S^1}(*)$.
From the discussion after Theorem A, we know that on the stalks
at zero $f^E_!(1)$ is the $S^1$-equivariant elliptic genus of
$X$, which is a priori a power series in $u$.  Theorem C with
$V=TX$ says that $1$ extends to a global section in
$\E^*_{S^1}(X)^{[f]}=\E^*_{S^1}(X)^{[TX]}$.  Therefore
$f^E_!(1)$ is the germ of a global section in 
$\E^*_{S^1}(*)=\cO_{\cE}$.  But any such section is a constant,
so the $S^1$-equivariant elliptic genus of $X$ is a constant.
This proves the rigidity of the elliptic genus
(Corollary~\ref{rigidity-theorem}).

Now the greater level of generality of Theorem C allows us
to extend the rigidity theorem to families of elliptic genera.
The question of stating and proving such a theorem was posed by
H. Miller in~\cite{Barcelona}.

\begin{thmd} \label{D}
(Rigidity for families) \textit{Let $\pi: E\ra B$ be a spin
oriented $S^1$-equivariant fibration.  Then the elliptic genus
of the family $\pi^E_!(1)$ is constant as a rational
function, i.e.\ when the generator $u$ of $H^*_{S^1}(*)=\C[u]$
is inverted.}
\end{thmd}

\newpage

\vs

\section{$S^1$-equivariant elliptic cohomology}
\label{section-ell-defn}

In this section we give the construction of $S^1$-equivariant
elliptic cohomology with complex coefficients.  But in order
to set up this functor, we need a few definitions.

\vs
\bc
{\sc 3.1. Definitions.}
\ec
\vs

Let $\cE$ be an elliptic curve over $\C$ with structure sheaf
$\cO_\cE$.  Let $\tt$ be a uniformizer of $\cE$, i.e.\
a generator of the maximal ideal of the local ring at zero
$\cO_{\cE,0}$. 
We say that $\tt$ is an \dfn{additive uniformizer} if
for all $x,y\in V_\tt$ such that $x+y\in V_\tt$, we have
$\tt(x+y)=\tt(x)+\tt(y)$.  An additive uniformizer always
exists, because we can take for example $\tt$ to be the
local inverse of the group map $\C \ra \cE$, where the
universal cover of $\cE$ is identified with $\C$.  Notice that
any two additive uniformizers differ by a nonzero constant,
because the only additive continuous functions on $\C$
are multiplications by a constant.

Let $V_\tt$ be a neighborhood of zero in $\cE$
such that $\tt:V_\tt \ra \C$ is a homeomorphism on its image.
Denote by $t_\aa$ translation by $\aa$ on $\cE$.  We say that
a neighborhood $V$ of $\aa\in\cE$ is \dfn{small} if 
$t_{-\aa}(V)\subseteq V_\tt$.

Let $\aa \in \cE$.
We say that $\aa$ is a torsion point of $\cE$ if there exists
$n>0$ such that $n\aa=0$.  The smallest $n$ with this property is
called the \dfn{exact order} of $\aa$.

Let $X$ be a finite $S^1$-$CW$ complex.
If $H \subseteq S^1$ is a subgroup, denote by $X^H$ the
submanifold of $X$ fixed by each element of $H$.
Let $\Z_n \subseteq S^1$ be the cyclic subgroup of order $n$.
Define a subgroup $H(\aa)$ of $S^1$ by: $H(\aa)=\Z_n$ if $\aa$ is
a torsion point of exact order $n$; $H(\aa)=S^1$ otherwise.
Then denote by 
   $$X^\aa=X^{H(\aa)}\period$$

Now suppose we are given an $S^1$-equivariant map
of $S^1$-$CW$ complexes $f:X\ra Y$.  A point $\aa \in \cE$ is called
\dfn{special} with respect to $f$ if either $X^\aa\neq X^{S^1}$ or
$Y^\aa\neq Y^{S^1}$.  When it is clear what $f$ is, we simply call
$\aa$ special.  A point $\aa\in\cE$ is called special with respect
to $X$ if it is special with respect to the identity function
$id:X\ra X$.

An indexed open cover $\cU=(U_\aa)_{\aa \in \cE}$ of $\cE$ is said to 
be \dfn{adapted} (with respect to $f$) if it satisfies the
following conditions:
\bi
\item[1.] $U_\aa$ is a small open neighborhood of $\aa$;
\item[2.] If $\aa$ is not special, then $U_\aa$ contains no special
   point;
\item[3.] If $\aa \neq \aa'$ are special points, 
   $U_\aa \cap U_{\aa'} = \emptyset$.
\ei
Notice that, if $X$ and $Y$ are finite $S^1$-$CW$ complexes,
then there exists an open cover of $\cE$ which is adapted to $f$.
Indeed, the set of special points is a finite subset of $\cE$.

If $X$ is a finite $S^1$-$CW$ complex, we define the
holomorphic $S^1$-equivariant cohomology of $X$ to be
  $$ \ho^*_{S^1}(X) =  H^*_{S^1}(X)\ts_{\C[u]}\cO_{\C,0} \period$$
$\cO_{\C,0}$ is the ring of germs of holomorphic functions at zero
in the variable $u$, or alternatively it is the subring of
$\C[\![u]\!]$ of convergent power series with positive radius of
convergence.

Notice that $\ho^*_{S^1}$ is not $\Z$-graded anymore, because
we tensored with the inhomogenous object $\cO_{\C,0}$.  However,
it is $\Z_2$-graded, by the even and odd part, because 
$\C[u]$ and $\cO_{\C,0}$ are concentrated in even degrees.

\vs
\bc
{\sc 3.2. Construction of} $\E^*_{S^1}$
\ec
\vs

We are going to define now a sheaf $\cF=\cF_{\tt,\cU}$ over $\cE$
whose stalk at $\aa \in \cE$ is isomorphic to $\ho^*_{S^1}(X^\aa)$.
Recall that, in order to give a sheaf
$\cF$ over a topological space, it is enough to give an open cover
$(U_\aa)_\aa$ of that space, and a sheaf $\cF_\aa$ on each $U_\aa$
together with isomorphisms of sheaves 
$\phi_{\aa\bb}: \cF_{\aa\mid_{U_\aa \cap U_\bb}} \lra 
                        \cF_{\bb\mid_{U_\aa \cap U_\bb}}$,
such that $\phi_{\aa\aa}$ is the identity function, and the cocycle
condition $\phi_{\bb\gg}\phi_{\aa\bb}=\phi_{\aa\gg}$ is satisfied on
$U_\aa \cap U_\bb \cap U_\gg$.

Fix $\tt$ an additive uniformizer of $\cE$.
Consider an adapted open cover $\cU=(U_\aa)_{\aa \in \cE}$.

\bdefn \label{local-sheaf}
Define a sheaf $\cF_\aa$ on $U_\aa$ by declaring for any open
$U \subseteq U_\aa$ 
  $$\cF_\aa(U) := 
       H^*_{S^1}(X^\aa)\ts_{\C[u]}\cO_{\cE}(U-\aa) \period$$
\edefn
\noindent
The map $\C[u] \ra \cO_{\cE}(U-\aa)$ is given by 
sending $u$ to $\tt$ (the germ $\tt$ extends to $U-\aa$ because
$U_\aa$ is small).  $U-\aa$ represents the translation of $U$ by
$-\aa$, and $\cO_{\cE}(U-\aa)$ is the ring of holomorphic functions
on $U-\aa$.  The restriction maps of the sheaf are defined so that
they come from those of the sheaf $\cO_{\cE}$.

First we notice that we can make $\cF_\aa$ into a sheaf of
$\cO_{\cE \;\mid U_\aa}$-modules: if $U \subseteq U_\aa$, we want
an action of $f\in\cO_{\cE}(U)$ on $\cF_\aa(U)$.  The translation
map $t_\aa: U-\aa \ra U$, which takes $u$ to $u+\aa$ gives a 
translation $t_\aa^*:\cO_{\cE}(U) \ra \cO_{\cE}(U-\aa)$, which 
takes $f(u)$ to $f(u+\aa)$.  Then we take the result of the action
of $f\in\cO_{\cE}(U)$ on 
$\mu\ts g \in \cF_\aa(U) = 
        H^*_{S^1}(X^\aa)\ts_{\C[u]}\cO_{\cE}(U-\aa)$
to be $\mu\ts (t^*_\aa f\cdot g)$.  Moreover, $\cF_\aa$ is coherent
because $H^*_{S^1}(X^\aa)$ is a finitely generated $\C[u]$-module.

Now for the second part of the definition of $\cF$, we have to 
glue the different sheaves $\cF_\aa$ we have just constructed.
If $U_\aa \cap U_\bb \neq \emptyset$ we need to define an isomorphism
of sheaves
$\phi_{\aa\bb}: \cF_{\aa\mid_{U_\aa \cap U_\bb}} \lra
                        \cF_{\bb\mid_{U_\aa \cap U_\bb}}$
which satisfies the cocycle condition.
Recall that we started with an adapted open cover 
$(U_\aa)_{\aa \in \cE}$.
Because of the condition 3 in the definition of an adapted cover,
$\aa$ and $\bb$ cannot be both special, so we only have to define 
$\phi_{\aa\bb}$ when, say, $\bb$ is not special.  In that case 
$X^\bb = X^{S^1}$.  Consider an arbitrary open set
$U \subseteq U_\aa \cap U_\bb$.

\bdefn \label{gluing-map}
Define $\phi_{\aa\bb}$ as the composite of the following 
maps:
\begin{equation}
 \ba{lcl}
\cF_\aa(U)& = &  H^*_{S^1}(X^\aa)\ts_{\C[u]}\cO_{\cE}(U-\aa)          \\
          &\ra&  H^*_{S^1}(X^\bb)\ts_{\C[u]}\cO_{\cE}(U-\aa)          \\
          &\ra&  (H^*(X^\bb) \ts_\C \C[u]) \ts_{\C[u]}\cO_{\cE}(U-\aa)\\
          &\ra&  H^*(X^\bb) \ts_\C \cO_{\cE}(U-\aa)                   \\
          &\ra&  H^*(X^\bb) \ts_\C \cO_{\cE}(U-\bb)                   \\
          &\ra&  (H^*(X^\bb) \ts_\C \C[u]) \ts_{\C[u]}\cO_{\cE}(U-\bb)\\
          &\ra&  H^*_{S^1}(X^\bb)\ts_{\C[u]}\cO_{\cE}(U-\bb)          \\
          & = &  \cF_\bb(U)  \period
 \ea   
\tag{$*$}
\end{equation}
The map on the second row is the natural map $i^* \ts 1$, where
$i:X^\bb \ra X^\aa$ is the inclusion.
\edefn

\blemma
$\phi_{\aa\bb}$ is an isomorphism.
\elemma

\bp
The second and and the sixth maps are isomorphisms because
$X^\bb=X^{S^1}$, and therefore 
$H^*_{S^1}(X^\bb)\we H^*(X^\bb) \ts_\C \C[u]$.  The properties of the
tensor product imply that the third and the fifth maps are
isomorphisms.  The fourth map comes from translation by
$\bb-\aa$, so it is also an isomorphism.

Finally, the second map $i^* \ts 1$ is an isomorphism because
\bi
\item[a)] If $\aa$ is not special, then $X^\aa=X^{S^1}=X^\bb$,
so $i^* \ts 1$ is the identity.
\item[b)] If $\aa$ is special, then $X^\aa \neq X^\bb$.
However, we have  $(X^\aa)^{S^1} = X^{S^1} = X^\bb$.
Then we can use the Atiyah--Bott localization theorem in equivariant
cohomology from~\cite{AB}.  This says that 
$i^*:H^*_{S^1}(X^\aa) \ra H^*_{S^1}(X^\bb)$
is an isomorphism after inverting $u$.  So it is enough to show
that $\tt$ is invertible in $\cO_{\cE}(U-\aa)$, because this would
imply that $i^*$ becomes an isomorphism after tensoring with 
$\cO_{\cE}(U-\aa)$ over $\C[u]$.  Now, because $\aa$ is special,
the condition 2 in the definition of an adapted cover says that
$\aa \notin U_\bb$.  But $U \subseteq U_\aa \cap U_\bb$, so 
$\aa \notin U$, hence $0 \notin U-\aa$.  This is equivalent to 
$\tt$ being invertible in $\cO_{\cE}(U-\aa)$.
\ei
\ep

\bremark \label{glue} To simplify notation, we can describe 
$\phi_{\aa\bb}$ as the composite of the following two maps:
%  $$\xymatrix{ H^*_{S^1}(X^\aa)\ts_{\C[u]}\cO_{\cE}(U-\aa) 
%                                     \ar[d]^{i^*}\\  
%               H^*_{S^1}(X^\bb)\ts_{\C[u]}\cO_{\cE}(U-\aa)
%                              \ar[d]^{\mbox{$t^*_{\bb-\aa}$}}\\
%               H^*_{S^1}(X^\bb)\ts_{\C[u]}\cO_{\cE}(U-\bb) \period 
%             }    $$
 $$  H^*_{S^1}(X^\aa)\ts_{\C[u]}\cO_{\cE}(U-\aa) 
                       \llra{i^*}  
     H^*_{S^1}(X^\bb)\ts_{\C[u]}\cO_{\cE}(U-\aa)
                       \llra{t^*_{\bb-\aa}}
     H^*_{S^1}(X^\bb)\ts_{\C[u]}\cO_{\cE}(U-\bb) \period    $$
By the first map we really mean $i^*\ts 1$.  
The second map is not $1\ts t^*_{\bb-\aa}$,
because $t^*_{\bb-\aa}$ is not a map of $\C[u]$-modules.  
However, we use $t^*_{\bb-\aa}$ as a shorthand for the
corresponding composite map specified in $(*)$.  Note that
$\phi_{\alpha\beta}$ is linear over $\cO_\cE(U)$, so we get
a map of sheaves of $\Z_2$-graded $\cO_\cE(U)$-algebras.
\eremark

One checks easily now that $\phi_{\aa\bb}$ satisfies the cocycle 
condition: Suppose we have three open sets $U_\aa$, $U_\bb$ and $U_\gg$
such that $U_\aa \cap U_\bb \cap U_\gg \neq \emptyset$.  Because our
cover was chosen to be adapted, at least two out of the three spaces
$X^\aa$, $X^\bb$ and $X^\gg$ are equal to $X^{S^1}$.
Thus the cocycle condition reduces essentially to
$t^*_{\gg-\bb} t^*_{\bb-\aa} = t^*_{\gg-\aa}$, which is clearly true.

\bdefn \label{f-sheaf}
Let $\cU=(U_\aa)_{\aa\in\cE}$ be an adapted cover of $\cE$,
and $\tt$ an additive uniformizer.  We define a sheaf
$\cF=\cF_{\tt,\cU}$ on $\cE$ by gluing the sheaves $\cF_\aa$
from Definition~\ref{local-sheaf} via the gluing maps
$\phi_{\aa\bb}$ defined in~\ref{gluing-map}.
\edefn

One can check now easily that $\cF$ is a coherent analytic sheaf of
algebras.  

Notice that we can remove the dependence of $\cF$
on the adapted cover $\cU$ as
follows:  Let $\cU$ and $\cV$ be two covers adapted to $(X,A)$.
Then any common refinement $\cW$ is going to be adapted as well, and
the corresponding maps of sheaves
$\cF_{\tt,\cU} \ra \cF_{\tt,\cW} \la \cF_{\tt,\cV}$ are isomorphisms
on stalks, hence isomorphisms of sheaves.  Therefore we can omit the
subscript $\cU$, and write $\cF=\cF_\tt$.
Next we want to show that $\cF_\tt$ is independent of the choice of
the additive uniformizer $\tt$.

\bprop
If $\tt$ and $\tt'$ are two additive uniformizers, then there exists
an isomorphism of sheaves of $\cO_\cE$-algebras
$f_{\tt\tt'}:\cF_\tt \ra \cF_{\tt'}$.  If $\tt''$ is a third
additive uniformizer, then $f_{\tt'\tt''} f_{\tt\tt'} = \pm f_{\tt\tt''}$.
\eprop

\bp
We modify slightly the notations used in Definition~\ref{local-sheaf}
to indicate the dependence on $\tt$: $\cF^\tt_\aa(U) = 
H^*_{S^1}(X^\aa)\ts^\tt_{\C[u]}\cO_{\cE}(U-\aa)$.
Recall that $u$ is sent to $\tt$ via the algebra map
$\C[u]\ra\cO_{\cE}(U-\aa)$.  If $\tt'$ is another additive uniformizer,
we saw at the beggining of this Section that there
exists a nonzero constant $a$ in $\cO_{\cE,0}$ such that
$\tt=a \tt'$.  Choose a square root of $a$ and denote it by $a^{1/2}$.
Define a map  $f_{\tt\tt',\aa}:\cF^\tt_\aa(U) \ra \cF^{\tt'}_\aa(U)$
by $x\ts^\tt g \mapsto a^{|x|/2} x\ts^{\tt'}g$.  We have assumed that
$x$ is homogeneous in $H^*_{S^1}(X^\aa)$, and that
$|x|$ is the homogeneous degree of $x$.

One can easily check that $f_{\tt\tt',\aa}$ is a map of sheaves
of $\cO_\cE$-algebras.  We also have
$\phi^{\tt'}_{\aa\bb} \circ f_{\tt\tt',\aa} = 
f_{\tt\tt',\bb} \circ \phi^\tt_{\aa\bb}$, which means that the maps
$f_{\tt\tt',\aa}$ glue to define a map of sheaves 
$f_{\tt\tt'}:\cF_\tt \ra \cF_{\tt'}$.
The equality $f_{\tt'\tt''} f_{\tt\tt'} = \pm f_{\tt\tt''}$ comes from
$(\tt'/\tt'')^{1/2} (\tt/\tt')^{1/2} = \pm (\tt/\tt'')^{1/2}$.
\ep

\bdefn
The $S^1$-equivariant elliptic cohomology of the finite $S^1$-$CW$
complex $X$ is the sheaf $\cF=\cF_{\tt,\cU}$ constructed above,
which according to the previous results does not depend
on the adapted open cover $\cU$ or on the
additive uniformizer $\tt$.  Denote this sheaf by $\E^*_{S^1}(X)$.
\edefn

If $X$ is a point, one can see that $\E^*_{S^1}(X)$ is the
structure sheaf $\cO_\cE$.

\bthm \label{cohomology-theory}
$\E^*_{S^1}(\blank)$ defines an $S^1$-equivariant cohomology theory
with values in the category of coherent analytic sheaves of
$\Z_2$-graded $\cO_\cE$-algebras.
\ethm

\bp
For $\E^*_{S^1}(\blank)$ to be a cohomology theory, we need
naturality.  Let $f:X\ra Y$ be an $S^1$-equivariant map of
finite $S^1$-$CW$ complexes.  We want to define a map of
sheaves $f^*:\E^*_{S^1}(Y)\ra \E^*_{S^1}(X)$ with the properties
that $1_X^* = 1_{\E^*_{S^1}(X)}$ and $(fg)^*=g^* f^*$.
Choose $\cU$ an open cover adapted to $f$, and $\tt$ an additive
uniformizer of $\cE$.  Since $f$ is $S^1$-equivariant, for each $\aa$
we get by restriction a map $f_\aa : X^\aa \ra Y^\aa$.
This induces a map $H^*_{S^1}(Y^\aa)\ts_{\C[u]}\cO_{\cE}(U-\aa)
\llra{f^*_\aa \ts 1} H^*_{S^1}(X^\aa)\ts_{\C[u]}\cO_{\cE}(U-\aa)$.
To get our global map $f^*$, we only have to check that
$f^*_\aa$ glue well, i.e.\ that they commute with the gluing maps
$\phi_{\aa\bb}$.  This follows easily from the naturality of ordinary
equivariant cohomology, and from the naturality in $X$ of the
isomorphism $H^*_{S^1}(X^{S^1}) \cong H^*(X^{S^1})\ts_\C \C[u]$.

Also, we need to define $\E^*_{S^1}$ for pairs.  Let $(X,A)$ be a pair
of finite $S^1$-$CW$ complexes, i.e.\ $A$ is a closed subspace of $X$,
and the inclusion map $A\ra X$ is $S^1$-equivariant.  We then define
$\E^*_{S^1}(X,A)$ as the kernel of the map
$j^*: \E^*_{S^1}(X/A) \ra \E^*_{S^1}(*)$, where $j:*=A/A \ra X/A$
is the inclusion map.  If $f:(X,A)\ra (Y,B)$ is a map of pairs of finite
$S^1$-$CW$ complexes, then $f^*:\E^*_{S^1}(Y,B)\ra \E^*_{S^1}(X,A)$ is
defined as the unique map induced on the corresponding kernels from
$f^*:\E^*_{S^1}(Y)\ra \E^*_{S^1}(X)$.

Now we have to define the coboundary map
$\dd:\E^*_{S^1}(A)\ra \E^{*+1}_{S^1}(X,A)$.  This is obtained by
gluing the maps $H^*_{S^1}(A^\aa)\ts_{\C[u]}\cO_{\cE}(U-\aa)
\llra{\dd_\aa \ts 1} H^{*+1}_{S^1}(X^\aa,A^\aa)\ts_{\C[u]}\cO_{\cE}(U-\aa)$,
where $\dd_\aa:H^*_{S^1}(A^\aa) \ra H^{*+1}_{S^1}(X^\aa,A^\aa)$ is
the usual coboundary map.  The maps $\dd_\aa\ts 1$ glue well, because
$\dd_\aa$ is natural.

To check the usual axioms of a cohomology theory: naturality, exact
sequence of a pair, and excision for $\E^*_{S^1}(\blank)$, recall
that this sheaf was obtained by gluing the sheaves $\cF_\aa$ along the
maps $\phi_{\aa\bb}$.  Since $\cF_\aa$ were defined using
$H^*_{S^1}(X^\aa)$, the properties of ordinary $S^1$-equivariant
cohomology pass on to $\E^*_{S^1}(\blank)$, as long as tensoring
with $\cO_{\cE}(U-\aa)$ over $\C[u]$ preserves exactness.  But this
is a classical fact: see for example the appendix of Serre~\cite{Ser}.
\ep

This proves THEOREM A stated in Section~\ref{section-results}.

\bremark
Notice that we can arrange our functor $\E^*_{S^1}(\blank)$ to take
values in the category of coherent \emph{algebraic} sheaves
over $\cE$ rather than in the category of coherent \emph{analytic} sheaves.  
This follows from a theorem of Serre~\cite{Ser} which says that the
the categories of coherent holomorphic sheaves and coherent algebraic
sheaves over a projective variety are equivalent.
\eremark

\vs
\bc
{\sc 3.3. Alternative description of} $\E^*_{S^1}$
\ec
\vs

For calculations with $\E^*_{S^1}(\blank)$ we want a
description which involves a finite open cover of $\cE$.
Start with an adapted open cover $(U_\aa)_{\aa \in \cE}$.
Recall that the set of special points with respect to $X$ is 
finite.  Denote this set by $\{\aa_1,\ldots,\aa_n\}$.
To simplify notation, denote for $i=1,\ldots,n$
  $$ U_i:=U_{\aa_i}, \mbox{ and }
     U_0:=\cE\setminus \{\aa_1,\ldots,\aa_n\} \period $$

On each $U_i$, with $0\leq i \leq n$, we define a sheaf $\cG$
as follows:
\bi
\item[a)] If $1 \leq i \leq n$, then $\forall U\subseteq U_i$,
 $\cG_i(U):=  H^*_{S^1}(X^{\aa_i})\ts_{\C[u]}\cO_{\cE}(U-\aa_i)$.
 The map $\C[u] \ra \cO_{\cE}(U-\aa_i)$ was described 
 in Definition~\ref{local-sheaf}.
\item[b)] If $i=0$, then $\forall U\subseteq U_0$,
 $\cG_i(U):= H^*(X^{S^1})\ts_\C \cO_{\cE}(U)$.
\ei
Now glue each $\cG_i$ to $\cG_0$ via the map of sheaves $\phi_{i0}$ 
defined as the composite of the following isomomorphisms 
($U\subseteq U_i\cap U_0$):
$H^*_{S^1}(X^{\aa_i})\ts_{\C[u]}\cO_{\cE}(U-\aa_i) 
                                    \llra{i^*\ts 1}              
              H^*_{S^1}(X^{S^1}) \ts_{\C[u]} \cO_{\cE}(U-\aa_i)
                             \llra{\cong} 
              H^*(X^{S^1}) \ts_\C \cO_{\cE}(U-\aa_i)
                             \llra{t^*_{-\aa_i}}     
              H^*(X^{S^1}) \ts_\C \cO_{\cE}(U)$.

Since there cannot be three distinct $U_i$ with nonempty intersection,
there is no cocycle condition to verify.

\bprop
The sheaf $\cG$ we have just described is isomorphic to $\cF$,
thus allowing an alternative definition of $\E^*_{S^1}(X)$.
\eprop

\bp
One notices that $U_0=\cup\{U_\bb\; |\; \bb \mbox{ nonspecial}\} $,
because of the third condition in the definition of an adapted cover.
If $U\subseteq \cup_\bb U_\bb$, a global section in $\cF(U)$
is a collection of sections $s_\bb \in \cF(U\cap U_\bb - \bb)$
which glue, i.e.\ $t^*_{\bb-\bb'}s_\bb = s_{\bb'}$.
So $t^*_{-\bb}s_\bb = t^*_{-\bb'}s_{\bb'}$ in 
$\cG(U\cap U_\bb \cap U_{\bb'})$, which means that we get an element
in $\cG(U)$, since the $U_\bb$'s cover $U$.  So 
$\cF_{\mid U_0} \cong \cG_{\mid U_0}$.  But clearly 
$\cF_{\mid U_i} \cong \cG_{\mid U_i}$ for $1 \leq i \leq n$, and the
gluing maps are compatible.  Therefore $\cF \cong \cG$.
\ep

As it is the case with any coherent sheaf of 
$\cO_\cE$-modules over an elliptic curve, 
$\E^*_{S^1}(X)$ splits (noncanonically) into a direct 
sum of a locally free sheaf, i.e.\ the sheaf of sections of some
holomorphic vector bundle, and a sum of skyscraper sheaves.

Given a particular $X$, we can be more specific:  We know that 
$H^{*}_{S^1}(X)$ splits noncanonically into a free
and a torsion $\C[u]$-module.  Given such a splitting, we can speak 
of the free part of $H^{*}_{S^1}(X)$.  Denote it by 
$H^{*}_{S^1}(X)_{free}$.  The map 
$i^* H^{*}_{S^1}(X)_{free} \ra H^{*}_{S^1}(X^{S^1})$
is an injection of finitely generated free $\C[u]$-modules of 
the same rank, by the localization theorem.
$\C[u]$ is a p.i.d., so by choosing appropriate bases in 
$H^{*}_{S^1}(X)_{free}$ and $H^{*}_{S^1}(X^{S^1})$, the map
$i^*$ can be written as a diagonal matrix 
$D(u^{n_1},\ldots,u^{n_k})$, $n_i\geq 0$.  Since $i^*1=1$, we can 
choose $n_1=0$.

So at the special points $\aa_i$, we have the map
$i^*:H^{*}_{S^1}(X^{\aa_i})_{free} \ra H^{*}_{S^1}(X^{S^1})$,
which in appropriate bases can be written as a diagonal matrix
$D(1,u^{n_2},\ldots,u^{n_k})$.  This gives over $U_i\cap U_0$ the
transition functions $u\mapsto D(1,u^{n_2},\ldots,u^{n_k})\in GL(n,\C)$.
However, we have to be careful since the basis of $H^{*}_{S^1}(X^{S^1})$
changes with each $\aa_i$, which means that the transition functions
are diagonal only up to a (change of base) matrix.  But this matrix
is invertible over $\C[u]$, so we get that the free part of
$\E^*_{S^1}(X)$ is a sheaf of sections of a holomorphic vector bundle.

An interesting question is what holomorphic vector bundles one
gets if $X$ varies.  Recall that holomorphic vector bundles over
elliptic curves were classified by Atiyah in 1957.

\bexample \label{example-thom-sheaf}
Calculate $\E^*_{S^1}(X)$ for $X=S^2(n)=$ the 2-sphere with the 
$S^1$-action which rotates $S^2$ $n$ times around the north-south 
axis as we go once around $S^1$.  If $\aa$ is an $n$-torsion
point, then $X^\aa=X$.  Otherwise, $X^\aa=X^{S^1}$,
which consists of two points: $\{P_+,P_-\}$, the North and the
South poles.  Now $H^*_{S^1}(S^2(n))=H^*(BS^1\vee BS^1)=
\C[u]\times_\C\C[u]$, on which $\C[u]$ acts diagonally.
$i^*: H^*_{S^1}(X) \ra H^*_{S^1}(X^{S^1})$ is the inclusion
$\C[u]\times_\C \C[u] \inc \C[u]\times\C[u]$.

Choose the bases 
\bi
\item[a)] $\{(1,1),(u,0)\}$ of $\C[u]\times_\C\C[u]$;
\item[b)] $\{(1,1),(1,0)\}$ of $\C[u]\times\C[u]$.
\ei
Then $H^*_{S^1}(X)\llra{\sim}\C[u]\oplus\C[u]$ by
$(P(u),Q(u))\mapsto (P,\frac{Q-P}u)$, and 
$H^*_{S^1}(X^{S^1})\llra{\sim}\C[u]\oplus\C[u]$ by
$(P(u),Q(u))\mapsto (P,Q-P)$.  Hence $i^*$ is given by the
diagonal matrix $D(1,u)$.  So $\E^*_{S^1}(X)$ looks locally
like $ \cO_{\C P^1} \oplus \cO_{\C P^1}(-1\cdot 0)$.  This happens 
at all the $n$-torsion points of $\cE$, so
$\E^*_{S^1}(X)\cong \cO_\cE \oplus \cO_\cE(\Delta)$, where
$\Delta$ is the divisor which consists of all $n$-torsion 
points of $\cE$, with multiplicity $1$.  

One can also check
that the sum of all $n$-torsion points is zero, so by
Abel's theorem the
divisor $\Delta$ is linearly equivalent to $-n^2\cdot 0$.
Thus $\E^*_{S^1}(S^2(n))\cong \cO_\cE \oplus \cO_\cE(-n^2\cdot 0)$.
We stress that the decomposition is only true as 
sheaves of $\cO_\cE$-modules, not as sheaves of $\cO_\cE$-algebras. 
\eexample

\bremark
Notice that $S^2(n)$ is the Thom space of the $S^1$-vector space
$\C(n)$, where $z$ acts on $\C$ by complex multiplication with $z^n$.
This means that the Thom isomorphism doesn't hold in $S^1$-equivariant
elliptic cohomology, because $\E^*_{S^1}(*)=\cO_\cE$, while
the reduced $S^1$-equivariant elliptic cohomology of the Thom space is
$\tilde{\E}^*_{S^1}(S^2(n)) = \cO_\cE(-n^2\cdot 0)$.
\eremark

\vs

\section{$S^1$-equivariant elliptic pushforwards}
\label{section-pushforwards}

While the construction of $\E^*_{S^1}(X)$ depends only on 
the elliptic curve $\cE$, the construction of the elliptic
pushforward $f^E_!$ involves extra choices, namely that of a
torsion point of exact order two on $\cE$, and a trivialization
of the cotangent space of $\cE$ at zero.

\vs
\bc
{\sc 4.1. The Jacobi sine}
\ec
\vs

Let $(\cE,P,\mu)$ be a triple formed with a nonsingular
elliptic curve $\cE$ over $\C$, a torsion point $P$ on $\cE$
of exact order two, and a 1-form $\mu$ which generates the
cotangent space $T^*_0\cE$.  For example, we can take
$\cE = \C/\Ll$, with $\Ll=\Z \oo_1 + \Z \oo_2$
a lattice in $\C$, $P=\oo_1/2$, and $\mu=dz$ at zero,
where $z$ is the usual complex coordinate on $\C$.

As in Hirzerbruch, Berger and Jung (\cite{HBJ}, Section 2.2),
we can associate to this data a function $s(z)$ on $\C$
which is elliptic (doubly periodic) with respect to a sublattice
$\tilde{\Ll}$ of index 2 in $\Ll$, namely
$\tilde{\Ll}=\Z \oo_1 + 2\Z \oo_2$.  (This leads to a double
covering $\tilde{\cE} \ra \cE$, and $s$ can be regarded as
a rational function on the ``doubled'' elliptic curve
$\tilde{\cE}$.)  Indeed, we can define $s$ up to a constant
by defining its divisor to be
   $$D = (0) + (\oo_1/2) - (\oo_2) - (\oo_1/2+\oo_2) \period$$
Then we can make $s$ unique by requiring that $ds = dz$ at zero.
We call this $s$ the \dfn{Jacobi sine}.  It has the following
properties (see~\cite{HBJ}):

\bprop \label{Jacobi-sine} \hspace{1mm} 
\bi
\item[a)] $s(z)$ is odd, i.e.\ $s(-z)=-s(z)$.  Around zero, $s$
can be expanded as a power series $s(z)=z+a_3z^3+a_5z^5+\cdots$.
\item[b)] $s(z+\oo_1)=s(z);s(z+\oo_2)=-s(z)$.
\item[c)] $s(z+\oo_1/2)=a/s(z)$, $a\neq 0$ (this follows
by looking at the divisor of $s(z+\oo_1/2)$).
\ei
\eprop

We now show that the construction of $s$ is canonical, i.e.\
it does not depend on the identification $\cE\cong\C/\Ll$.

\bprop
The definition of $s$ only depends on the triple $(\cE,P,\mu)$.
\eprop

\bp
First, we show that the construction of 
$\tilde{\cE}=\C/\tilde{\Ll}$ is canonical:
Let $\cE\cong\C/\Ll'$ be another identification of $\cE$.
We then have $\Ll'=\Z\oo_1' + \Z\oo_2'$, and $P$ is identified
with $\oo_1'/2$.  Since $\cE$ is also identified with $\C/\Ll$,
we get a group map $\ll:\C/\Ll \we \C/\Ll'$.  This implies that we
have a continuous group map $\ll:\C \we \C$ such that $\ll(\Ll)=\Ll'$.
Any such map must be multiplication by a nonzero constant $\ll\in\C$.
Moreover, we know that $\ll\oo_1/2=\oo_1'/2$.  This implies
$\ll\oo_1=\oo_1'$, and since $\ll$ takes $\Ll$ isomorphically
onto $\Ll'$, it follows that $\ll\oo_2=\pm\oo_2' + m\oo_1'$
for some integer $m$.  Multiplying this by 2, we get
$\ll\cdot2\oo_2=\pm2\oo_2' + 2m\oo_1'$.  This, together with
$\ll\oo_1=\oo_1'$, imply that multiplication by $\ll$ descends
to a group map  $\C/\tilde{\Ll} \we \C/\tilde{\Ll'}$.
But this precisely means that the construction of 
$\tilde{\cE}$ is canonical.

Notice that $P$ can be thought canonically as a point on
the ``doubled'' ellptic curve $\tilde{\cE}$.  We denote
by $P_1$ and $P_2$ the other two points of exact order 2
on $\tilde{\cE}$.  Then we form the divisor
   $$D = (0) + (P) - (P_1) - (P_2) \period$$
Although the choice of $P_1$ and $P_2$ is
noncanonical, the divisor $D$ is canonical, i.e.\ 
depends only on $P$.  Let $s$ be an elliptic function
on $\tilde{\cE}$ associated to the divisor $D$.  The
choice of $s$ is well-defined up to a constant which can
be fixed if we require that $ds=\pi^*\mu$ at zero, where
$\pi:\tilde{\cE} \ra \cE$ is the projection map.
\ep

\vspace{5mm}

Next, we start the construction of $S^1$-equivariant
elliptic pushforwards.  Let $f: X \ra Y$ be an equivariant
map between compact $S^1$-manifolds such that the restrictions
$f: X^\aa \ra Y^\aa$ are oriented maps.  Then we follow
Grojnowski~\cite{Gr} and define the pushforward of $f$
to be a map of sheaves 
$f^E_!:\E^*_{S^1}(X)^{[f]} \ra \E^*_{S^1}(Y)$,
where $\E^*_{S^1}(X)^{[f]}$ is the sheaf $\E^*_{S^1}(X)$ 
twisted by a $1$-cocycle to be defined later.

The main technical ingredient in the construction of the 
(global i.e.\ sheafwise) elliptic pushforward 
$f^E_!:\E^*_{S^1}(X)^{[f]} \ra \E^*_{S^1}(Y) \comma$
is the (local i.e.\ stalkwise) elliptic pushforward
$f^E_!:\ho^*_{S^1}(X^\aa) \ra \ho^*_{S^1}(Y^\aa)$.

In the following subsection, we construct elliptic Thom
classes and elliptic pushforwards in $\ho^*_{S^1}(\blank)$.
The construction is standard, with the only problem that
in order to show that something belongs to 
$\ho^*_{S^1}(\blank)$, we need some holomorphicity results
on characteristic classes.

\vs
\bc
{\sc 4.2. Preliminaries on pushforwards}
\ec
\vs

Let $\pi:V\ra X$ be a $2n$-dimensional oriented real $S^1$-vector
bundle over a finite $S^1$-$CW$ complex $X$, i.e.\ a vector bundle
with a linear action of $S^1$, such that $\pi$ commutes with the
$S^1$ action.  Now, for any space $A$ with an $S^1$ action, we can
define its Borel construction $A\times_{S^1}ES^1$, where $ES^1$ is
the universal principal $S^1$-bundle.  This construction is
functorial, so we get a vector bundle $V_{S^1}$ over $X_{S^1}$.
This has a classifying map $f_V:X_{S^1} \ra BSO(2n)$.
If $V_{univ}$ is the universal orientable vector bundle over
$BSO(2n)$, we also have a map of pairs, also denoted by
$f_V:(DV_{S^1},SV_{S^1})\ra (DV_{univ},SV_{univ})$.  As usual,
$DV$ and $SV$ represent the disc and the sphere bundle of $V$, 
respectively.

But it is known that the pair $(DV_{univ},SV_{univ})$ is
homotopic to $(BSO(2n),BSO(2n-1))$.  Also, we know that 
$H^*BSO(2n)=\C[p_1,\ldots,p_n,e] / (e^2-p_n)$,
where $p_j$ is the universal $j$'th Pontrjagin class, and 
$e$ is the universal Euler class.  From the long exact 
sequence of the pair, it follows that
$H^*(BSO(2n),BSO(2n-1))$ can be regarded as the ideal 
generated by $e$ in $H^*BSO(2n)$.  
The class $e\in H^*(DV_{univ},SV_{univ})$ is the universal
Thom class, which we will denote by $\phi_{univ}$. 
Then the ordinary equivariant Thom class of $V$ is defined
as the pullback class $f_V^* \phi_{univ} \in H^*_{S^1}(DV,SV)$,
and we denote it by $\phi_{S^1}(V)$.
Denote by $H^{**}_{S^1}(X)$ the completion of the module
$H^*_{S^1}(X)$ with respect to the ideal generated by $u$
in $H^*(BS^1)=\C[u]$.

Consider the power series $Q(x)= s(x)/x$, where $s(x)$ is the Jacobi
sine.  Since $Q(x)$ is even, Definition~\ref{defn-real-equiv-char-class}
gives a class $\mu_Q(V)_{S^1} \in H^{**}_{S^1}(X)$.  Then we
define a class in $H^{**}_{S^1}(DV,SV)$ by
$\phi_{S^1}^E(V) = \mu_Q(V)_{S^1} \cdot \phi_{S^1}(V)$.
One can also say that $\phi^E_{S^1}(V)=s(x_1)\cdots s(x_n)$, 
while $\phi_{S^1}(V)=x_1\ldots x_n$, where 
$x_1, \ldots, x_n$ are the equivariant Chern roots of $V$.
We call $\phi^E_{S^1}(V)$ the
\dfn{elliptic equivariant Thom class} of $V$.

Also, we define $e^E_{S^1}(V)$, the
equivariant elliptic Euler class
of $V$, as the image of $\phi^E_{S^1}(V)$ via the restriction map
$H^{**}_{S^1}(DV,SV) \ra H^{**}_{S^1}(X)$.

\bprop
If $V\ra X$ is an even dimensional real oriented $S^1$-vector
bundle, and $X$ is a finite $S^1$-$CW$ complex, then 
$\phi^E_{S^1}(V)$ actually lies in $\ho^*_{S^1}(DV,SV)$.
Cup product with the elliptic Thom class
$$\xymatrix{
 \ho^*_{S^1}(X) \ar[rr]^{\cup\,\phi^E_{S^1}(V)} & & \ho^*_{S^1}(DV,SV) \comma
           }   $$
is an isomorphism, the Thom isomorphism in $\ho$-theory.
\eprop

\bp
The difficult part, namely that $\mu_Q(V)_{S^1}$ is holomorphic,
is proved in the Appendix, in Proposition~\ref{holomorphic}.
Consider the usual cup product, which is a map
$\cup:H^*_{S^1}(X) \ts H^*_{S^1}(DV,SV) \ra H^*_{S^1}(DV,SV)$, and extend
it by tensoring with $\cO_{\C,0}$ over $\C[u]$.  We obtain a map
$\cup:\ho^*_{S^1}(X) \ts \ho^*_{S^1}(DV,SV) \ra \ho^*_{S^1}(DV,SV)$.
The equivariant elliptic Thom class of $V$ is
$\phi^E_{S^1}(V) = \mu_Q(V)_{S^1} \cup \phi_{S^1}(V)$,
so we have to show that both these classes are holomorphic.
But by Proposition~~\ref{holomorphic} in the Appendix, 
$\mu_Q(V)_{S^1}\in \ho^*_{S^1}(X)$.  And the ordinary Thom class
$\phi_{S^1}(V)$ belongs to $H^*_{S^1}(DV,SV)$, so it also belongs
to the larger ring $\ho^*_{S^1}(DV,SV)$.

Now, cup product with $\phi^E_{S^1}(V)$ gives an isomorphism
because $Q(x)=s(x)/x$ is an invertible power series around zero.
\ep

\bcor \label{local-pushforward}
If $f:X \ra Y$ is an $S^1$-equivariant oriented map between
compact $S^1$-manifolds, then there is an elliptic pushforward
  $$ f^E_!:\ho^*_{S^1}(X) \ra \ho^*_{S^1}(Y) \comma$$
which is a map of $\ho^*_{S^1}(Y)$-modules.
In the case when $Y$ is a point, $f^E_!(1)$ is the $S^1$-equivariant
elliptic genus of $X$.
\ecor

\bp
Recall (Dyer~\cite{Dy}) that the ordinary pushforward is defined as 
the composition of three maps, two of which are Thom isomorphisms,
and the third is a natural one.  The existence of the elliptic
pushforward follows therefore from the previous corollary.  The proof
that $f^E_!$ is a map of $\ho^*_{S^1}(Y)$-modules is the same as for
the ordinary pushforward.

The last statement is an easy consequence of the topological
Riemann--Roch theorem (see again \cite{Dy}), and of the definition 
of the equivariant elliptic Thom class.
\ep

Notice that, if $Y$ is point, $\ho^*_{S^1}(Y) \cong \cO_{\C,0}$, so
the $S^1$-equivariant elliptic genus of $X$ is holomorphic around
zero.  Also, if we replace 
$\ho^*_{S^1}(\blank) = H^*_{S^1}(\blank)\ts_{\C[u]}\cO_{\C,0}$
by $\hm^*_{S^1}(\blank) = H^*_{S^1}(\blank)\ts_{\C[u]} \cM(\C)$, where
$\cM(\C)$ is the ring of global meromorphic functions on $\C$,
the same proof as above shows that the $S^1$-equivariant elliptic
genus of $X$ is meromorphic in $\C$.  This proves the following
result, which is PROPOSITION B stated in Section~\ref{section-results}.

\bprop
The $S^1$-equivariant elliptic genus of a compact
$S^1$-manifold is the Taylor expansion at zero of a function on
$\C$ which is holomorphic at zero and meromorphic everywhere.
\eprop

\bc
{\sc 4.3. Construction of} $f^E_!$
\ec
\vs

The local construction of elliptic pushforwards is completed.
We want now to assemble the pushforwards in a map of
sheaves.  Let $f:X\ra Y$ be a map of compact $S^1$-manifolds which
commutes with the $S^1$-action.  We assume that either $f$ is
complex oriented or spin oriented, i.e.\ that the stable normal
bundle in the sense of Dyer~\cite{Dy} is complex oriented or spin
oriented, respectively.  (Grojnowski treats only the complex
oriented case, but in order to understand rigidity we also
need the spin case.)

Let $\cU$ be an open cover of $\cE$ adapted to $f$.  Let
$\aa,\bb \in \cE$ be such that $U_\aa\cap U_\bb \neq \emptyset$.
This implies that at least one point, say $\bb$, is nonspecial,
so $X^\bb=X^{S^1}$ and
$Y^\bb=Y^{S^1}$.  We specify now the orientations of the maps
and vector bundles involved.  Since $X^\bb=X^{S^1}$, the normal
bundle of the embedding $X^\bb \ra X^\aa$ has a complex structure,
where all the weights of the $S^1$-action on $V$ are positive.

If $f$ is complex oriented, it follows that the restriction maps
$f^\aa:X^\aa \ra Y^\aa$ and $f^\bb:X^\bb \ra Y^\bb$ are also complex
oriented, hence oriented.
If $f$ is spin oriented, this means that the stable normal
bundle $W$ of $f$ is spin.  If $H$ is any subgroup of $S^1$, we know that
the vector bundle $W^H \ra X^H$ is oriented: If $H=S^1$, $W$ splits
as a direct sum of $W^H$ with a bundle corresponding to the nontrivial
irreducible representations of $S^1$; this latter bundle is
complex, hence oriented, so the orientation of $W$ induces one on
$W^H$.  If $H=\Z_n$, Lemma 10.3 of Bott and Taubes~\cite{BT} implies
that $W^H$ is oriented.  In conclusion, both maps $f^\aa$ and $f^\bb$
are oriented.

According to Corollary~\ref{local-pushforward}, we can define
elliptic pushforwards at the level of stalks:
$(f^\aa)^E_!:\ho^*_{S^1}(X^\aa) \ra \ho^*_{S^1}(Y^\aa)$
and $(f^\bb)^E_!:\ho^*_{S^1}(X^\bb) \ra \ho^*_{S^1}(Y^\bb)$.
The problem is that pushforwards do not commute with pullbacks, i.e.\
if $i:X^\bb\ra X^\aa$ and $j:Y^\bb\ra Y^\aa$ are the inclusions,
then it is not true in general that $j^*(f^\aa)^E_!=(f^\bb)^E_!i^*$.
However, by twisting the maps with some appropriate Euler classes,
the diagram becomes commutative.  Denote by 
$e^E_{S^1}(X^\aa/X^\bb)$ the $S^1$-equivariant Euler class of the
normal bundle to the embedding $i$, and by $e^E_{S^1}(Y^\aa/Y^\bb)$
the $S^1$-equivariant Euler class of the normal bundle to $j$.
Denote by 
 $$\ll_{\aa\bb}^{[f]}= 
    e^E_{S^1}(X^\aa/X^\bb)^{-1}\cdot (f^\bb)^* e^E_{S^1}(Y^\aa/Y^\bb)
           \period$$
A priori $\ll_{\aa\bb}^{[f]}$ belongs to the ring
$\ho^*_{S^1}(X^\bb)[\frac1{e^E_{S^1}(X^\aa/X^\bb)}]$, but we
will see later that we can improve this.

\blemma \label{commutativity-up-to-Euler-classes}
In the ring
$\ho^*_{S^1}(X^\bb)[\frac1u,\frac1{e^E_{S^1}(X^\aa/X^\bb)}]$
we have the following formula
  $$j^*(f^\aa)^E_!\mu^\aa = 
        (f^\bb)^E_!(i^*\mu^\aa\cdot \ll_{\aa\bb}^{[f]}) \comma $$
\elemma

\bp
From the hypothesis, we know that $i^*i^E_!$ is an isomorphism,
because it is multiplication by the invertible class 
$e^E_{S^1}(X^\alpha/X^\beta)$.  Also, since $u$ is invertible,
the localization theorem implies that $i^*$ is an isomorphism.
Therefore $i^E_!$ is an isomorphism.
Start with a class $\mu^\alpha$ on $X^\aa$.  Because $i^E_!$ is an
isomorphism, $\mu^\alpha$ can be written as $i^E_! \mu^\beta$, where
$\mu^\beta$ is a class on $X^\bb$.

Now look at the two sides of the equation to be proved:
\bi
\item[1.] The left hand side = $j^*(f^\aa)^E_!i^E_!\mu^\beta =
 j^*j^E_!(f^\bb)^E_!\mu^\beta =
 (f^\bb)^E_!\mu^\beta\cdot e^E_{S^1}(Y^\alpha/Y^\beta)$,
 because $j^*j^E_!=$ multiplication by $e^E(Y^\alpha/Y^\beta)$.
\item[2.] The right hand side = $(f^\bb)^E_![i^*i^E_!\mu^\beta\cdot 
 e^E_{S^1}(X^\alpha/X^\beta)^{-1}\cdot (f^\bb)^* e^E_{S^1}(Y^\alpha/Y^\beta)]=
 (f^\bb)^E_![\mu^\beta\cdot (f^\bb)^* e^E_{S^1}(Y^\alpha/Y^\beta)]=
 (f^\bb)^E_!\mu^\beta\cdot e^E_{S^1}(Y^\alpha/Y^\beta)$,
 where the last equality comes from the fact that $(f^\bb)^E_!$ is a map
 of $\ho^*_{S^1}(Y^\bb)$-modules.
\ei
\ep

Let $f:X\ra Y$ be a complex or spin oriented $S^1$-map.  Let $\cU$ 
be an open cover adapted to $f$, and $\aa,\bb\in\cE$ such that
$U_\aa\cap U_\bb \neq \emptyset$.  We know that $\aa$ and $\bb$ cannot
be both special, so assume $\bb$ nonspecial.
Let $U\subseteq U_\aa\cap U_\bb$.  Since $\cU$ is adapted, $\aa\notin U$.

\bprop  \label{commutative-lambda-diagram}
With these hypotheses,
$\ll_{\aa\bb}^{[f]}$ belongs to $H^*_{S^1}(X^\bb)\ts_{\C[u]}\cO_\cE(U-\bb)$,
and the following diagram is commutative:
 $$\xymatrix{ H^*_{S^1}(X^\aa)\ts_{\C[u]}\cO_\cE(U-\aa) \ar[r]^{(f^\aa)^E_!} 
                                \ar[d]_{\ll_{\aa\bb}^{[f]}\cdot i^*} & 
        H^*_{S^1}(Y^\aa)\ts_{\C[u]}\cO_\cE(U-\aa)       \ar[d]^{j^*}  \\
        H^*_{S^1}(X^\bb)\ts_{\C[u]}\cO_\cE(U-\aa)  
                   \ar[r]^{(f^\bb)^E_!} \ar[d]_{t^*_{\bb-\aa}}  & 
        H^*_{S^1}(Y^\bb)\ts_{\C[u]}\cO_\cE(U-\aa) \ar[d]^{t^*_{\bb-\aa}}\\
        H^*_{S^1}(X^\bb)\ts_{\C[u]}\cO_\cE(U-\bb)  \ar[r]^{(f^\bb)^E_!}   & 
        H^*_{S^1}(Y^\bb)\ts_{\C[u]}\cO_\cE(U-\bb)
     }$$
\eprop

\bp
Denote by $W$ the normal bundle of the embedding 
$X^\bb=X^{S^1}\ra X^\aa$.  Let us show that, if $\aa\notin U$, then
$e^E_{S^1}(W)$ is invertible in 
$H^*_{S^1}(X^\bb)\ts_{\C[u]}\cO_\cE(U-\aa)$.  Denote by $w_i$ the
nonequivariant Chern roots of $W$, and by $m_i$ the corresponding
rotation numbers of $W$ (see Proposition~\ref{formula-Chern-roots}
in the Appendix).
Since $X^\bb=X^{S^1}$, $m_i \neq 0$.  Also, the $S^1$-equivariant
Euler class of $W$ is given by
  $$e_{S^1}(W)=(w_1+m_1u)\ldots (w_r+m_ru)=
    m_1\ldots m_r(u+w_1/m_1)\ldots (u+w_r/m_r) \period$$
But $w_i$ are nilpotent, so $e_{S^1}(W)$ is invertible as long as
$u$ is invertible.  Now $\aa\notin U$ translates to $0\notin U-\aa$,
which implies that the image of $u$ via the map
$\C[u] \ra \cO_\cE(U-\aa)$ is indeed invertible.  To deduce now
that $e^E_{S^1}(W)$, the elliptic $S^1$-equivariant Euler class of $W$,
is also invertible, recall that $e^E_{S^1}(W)$ and $e_{S^1}(W)$
differ by a class defined using the power series
$s(x)/x=1+a_3x^2+a_5x^4+\cdots$, which is invertible for $U$ small
enough.

So $\ll_{\aa\bb}^{[f]}$ exists, and by the previous Lemma, the upper
part of our diagram is commutative.  The lower part is trivially
commutative.
\ep

Now, since $i^*$ are essentially the gluing maps in the sheaf
$\cF=\E^*_{S^1}(X)$, we think of the maps $\ll_{\aa\bb}^{[f]}\cdot i^*$
as giving the sheaf $\cF$ twisted by the cocycle $\ll_{\aa\bb}^{[f]}$.
Recall from Definition~\ref{f-sheaf} that $\cF$ was obtained
by gluing the sheaves $\cF_\aa$ over an 
adapted open cover $(U_\aa)_{\aa \in \cE}$.

\bdefn \label{twisted-glue}
The twisted gluing functions $\phi_{\aa\bb}^{[f]}$
are defined as the composition of the following three maps
$H^*_{S^1}(X^\aa)\ts_{\C[u]}\cO_\cE(U-\aa) \llra{i^*\ts 1}
H^*_{S^1}(X^\bb)\ts_{\C[u]}\cO_\cE(U-\aa) \llra{\cdot  \, \ll_{\aa\bb}^{[f]}}
H^*_{S^1}(X^\bb)\ts_{\C[u]}\cO_\cE(U-\bb) \llra{t^*_{\bb-\aa}}
H^*_{S^1}(X^\bb)\ts_{\C[u]}\cO_\cE(U-\bb)$.
The third map is defined as in Remark~\ref{glue}.
\edefn

As in the discussion after Remark~\ref{glue}, we can show
easily that $\phi_{\aa\bb}^{[f]}$ satisfy the cocycle condition.

\bdefn \label{pushforward}
Let $f:X\ra Y$ be an equivariant map of compact $S^1$-manifolds,
such that it is either complex or spin oriented.
We denote by $\E^*_{S^1}(X)^{[f]}$ the sheaf obtained
by gluing the sheaves $\cF_\aa$ defined in~\ref{local-sheaf},
using the twisted gluing functions $\phi_{\aa\bb}^{[f]}$.

Also, we define the $S^1$-equivariant elliptic pushforward of $f$
to be the map of coherent sheaves over $\cE$
   $$ f^E_!:\E^*_{S^1}(X)^{[f]} \ra \E^*_{S^1}(Y) $$
which comes from gluing the local elliptic pushforwards
$(f^\aa)^E_!$ (as defined in~\ref{local-pushforward}). 
We call $f^E_!$ the \dfn{Grojnowski pushforward}.
\edefn

The fact that $(f^\aa)^E_!$ glue well comes from the
commutativity of the diagram in
Proposition~\ref{commutative-lambda-diagram}.
The Grojnowski pushforward is functorial:
see~\cite{GKV} and~\cite{Gr}.

\vs

\section{Rigidity of the elliptic genus} \label{section-rigidity}

In this section we discuss the rigidity phenomenon in the context
of equivariant elliptic cohomology.  We start with a discussion
about orientations.

\vs
\bc
{\sc 5.1. Preliminaries on orientations}
\ec
\vs

Let $V\ra X$ be an even dimensional spin $S^1$-vector bundle over
a finite $S^1$-$CW$ complex $X$ (which means that the $S^1$-action
preserves the spin structure).  Let $n\in\N$.  We think of 
$\Z_n\subset S^1$ as the ring of $n$'th roots of unity in $\C$.
The invariants of $V$ under the actions of $S^1$ and $\Z_n$
are the $S^1$-vector bundles $V^{S^1} \ra X^{S^1}$ and
$V^{\Z_n}\ra X^{\Z_n}$.  We have $X^{S^1}\subseteq X^{\Z_n}$.

Let $N$ be a connected component of $X^{S^1}$, and $P$
a connected component of $X^{\Z_n}$ which contains $N$.
From now on we think of $V^{S^1}$ as a bundle over $N$,
and $V^{\Z_n}$ as a bundle over $P$.

Define the vector bundles $V/V^{S^1}$ and $V^{\Z_n}/V^{S^1}$
over $N$ by
 $$ V_{\mid N} = V^{S^1}\oplus V/V^{S^1} ; \;
    V^{\Z_n}_{\mid N}= V^{S^1}\oplus V^{\Z_n}/V^{S^1} \period$$
The decompositions of these two bundles come from the 
fact that $S^1$ acts trivially on the base $N$, so fibers
decompose into a trivial and nontrivial part.

Similarly, the action of $\Z_n$ on $P$ is trivial, so we get
a fiberwise decomposition of $V_{\mid P}$ by the different
representations of $\Z_n$:
    $$V_{\mid P} = 
              V^{\Z_n} \oplus 
                \bigoplus_{0<k<\frac{n}{2}}V(k) \oplus 
                   V(\frac{n}{2}) \period$$
By convention, $V(\frac n2)=0$ if $n$ is odd.
Lemma 9.4 in Bott and Taubes~\cite{BT} implies that 
$V^{\Z_n}$ and $V(\frac{n}2)$ are even dimensional
real oriented vector bundles.  Denote by
    $$V(K)=\bigoplus_{0<k<\frac{n}{2}}V(k)\period $$
Then we have the following decompositions:

\begin{equation} \label{VP}
   V_{\mid P} = V^{\Z_n} \oplus V(K) \oplus V(\frac{n}{2}) \period
\end{equation}

\begin{equation} \label{VN}
   V^{\Z_n}_{\mid N} = V^{S^1} \oplus V^{\Z_n}/V^{S^1}  \period
\end{equation}

Now we define the orientations for the different bundles involved:

First, if a bundle is oriented, any restriction to a smaller
base gets an induced orientation.

$V$ is oriented by its spin structure.

$\Z_n$ preserves the spin structure of $V$, so we can apply
Lemma 10.3 from Bott and Taubes~\cite{BT}, and deduce that
$V^{\Z_n}$ has an induced orientation.

Each $V(k)$ for $0<k<\frac{n}2$ has a complex structure,
for which  $g=e^{2\pi i/n} \in \Z_n$ acts by complex
multiplication with $g^k$.  So $V(K)$ has a complex orientation,
too.  Define the orientation on $V(K)$ by:
\bi
\item[$\bullet$] If $V(\frac n2)\neq 0$, $V(K)$ is oriented by 
 its complex structure described above. 
\item[$\bullet$] If $V(\frac n2)= 0$, then choose the orientation
 on $V(K)$ induced by the decomposition in~(\ref{VP}):
 $V_{\mid P} = V^{\Z_n} \oplus V(K)$.
\ei
Then the decomposition in~(\ref{VP}) induces an orientation
on $V(\frac n2)$.

We now orient bundles appearing in~(\ref{VN}) as follows:
Notice that $V^{\Z_n}/V^{S^1}$ has only nonzero rotation numbers,
so it has a complex structure for which all rotation numbers
are positive.
Define the orientation on $V^{\Z_n}/V^{S^1}$ by:
\bi
\item[$\bullet$] If $V^{S^1}\neq 0$, $V^{\Z_n}/V^{S^1}$ is
 oriented by its complex structure described above.
\item[$\bullet$] If $V^{S^1} = 0$, then 
 $V^{\Z_n}/V^{S^1} = V^{\Z_n}_{\mid N}$, so we choose 
 this orientation, induced from that on $V^{\Z_n}$ described
 above.
\ei

Finally, we orient $V/V^{S^1}$ from the decomposition
\begin{equation} \label{VS}
  V/V^{S^1} = V^{\Z_n}/V^{S^1} \oplus 
               V(K)_{\mid N} \oplus
                V(\frac{n}{2})_{\mid N} \period
\end{equation}

As a notational rule, we are going to use the subscript ``$or$'' to 
indicate the ``correct'' orientation on the given vector space, i.e.\
the orientations which we defined above.
When we omit the subscript ``$or$'', we assume the bundle has the
correct orientation.  But all bundles that appear in~(\ref{VS})
also have a complex structure (they have nonzero rotation numbers).
The subscript ``$cx$'' will indicate that we chose a
complex structure on the given vector space.  This is only
intended to make calculations easier.  So we choose complex
structures as follows:

For $V^{\Z_n}/V^{S^1}$ we choose as above the complex
structure where all rotation numbers are positive,
and similarly for $V(\frac{n}{2})_{\mid N}$.
Also, $V(K)_{\mid N}$ gets an induced complex structure
from the complex structure on $V(K)$ described above.
Now $V/V^{S^1}$ gets its complex structure from the 
decomposition~(\ref{VS}).

Let $i:N\ra P$ be the inclusion.  Here is a table with
the bundles of interest: 

\vskip5mm
\bc
\bt{|c|c|c|} 
\hline
bundle with the         & bundle with the       & sign difference between\\
correct orientation     & complex orientation   & the two orientations  \\
\hline 
$(V/V^{S^1})_{or}$   & $(V/V^{S^1})_{cx}$ &   $(-1)^{\ss}$       \\
$(V^{\Z_n}/V^{S^1})_{or}$ & $(V^{\Z_n}/V^{S^1})_{cx}$ 
                                  & $(-1)^{\ss(0)}$                 \\
$V(K)_{or}$          & $V(K)_{cx}$        & $(-1)^{\ss(K)}$        \\
$i^*(V(\frac{n}{2})_{or})$ & $(i^*V(\frac{n}{2}))_{cx}$ 
                                  & $(-1)^{\ss(\frac{n}{2})}$  \\
\hline
\et     
\ec
\vspace{4mm}

From the decomposition in (\ref{VS}) under the correct and the complex
orientations, we deduce that
\begin{equation} \label{signs}
  (-1)^{\ss(0)}(-1)^{\ss(K)} (-1)^{\ss(\frac{n}{2})} = (-1)^\ss \period 
\end{equation}

By the splitting principle (Bott and Tu~\cite{BTu}),
the pullback of $V/V^{S^1}$ to the flag manifold can be thought of as
a sum of complex line bundles $L(m_j)$, $j=1,\ldots,r$.  The complex
structure of $L(m_j)$ is such that $g\in S^1$ acts on $L(m_j)$ via
complex multiplication with $g^{m_j}$.  The numbers $m_j\neq 0$,
$j=1,\ldots,r$, are the rotation numbers.
By the real splitting principle, they are defined
also for even dimensional real oriented vector bundles, but in that
case the $m_j$'s are well defined only up to an even number of sign
changes.  We choose two systems of rotation numbers for 
$V/V^{S^1}$, one denoted by $(m_j)_j$, corresponding to
$(V/V^{S^1})_{or}$; and one denoted by $(m_j^*)_j$, corresponding to
$(V/V^{S^1})_{cx}$.  Of course, since the two orientations
differ by the sign $(-1)^\ss$, the systems $(m_j)_j$ and
$(m_j^*)_j$ will be the same up to a permutation
and a number of sign changes of the same parity with $(-1)^\ss$.

For $j=1,\ldots,r$, we define $q_j$ and $r_j$ as the quotient
and the remainder, respectively, of $m_j$ modulo $n$.  Similarly,
$q_j^*$ and $r_j^*$ are the quotient and the reminder of
$m_j^*$ modulo $n$.

We define now for each $k$ a set of indices of the corresponding
rotation numbers from the decomposition in~(\ref{VS}):
if $0\leq k \leq \frac{n}{2}$, define
   $$ I_k = \{j \in 1,\ldots,r\ |\ r_j= k \mbox{ or } n-k \}\period$$
Notice that for $0<k\leq \frac{n}2$, $I_k$ contains exactly the
indices of the rotation numbers for $V(k)$, and for $k=0$,
$I_0$ contains the indices of the rotation
numbers corresponding to $V^{\Z_n}/V^{S^1}$.
Also, define
  $$ I_K = \bigcup_{0<k<\frac{n}2} I_k  \period$$

% ** This is an attempt for an easy proof **
% By the splitting principle (Bott and Tu~\cite{BTu}),
% the pullback of $V/V^{S^1}$ to the flag manifold can be thought of as
% a sum of complex line bundles $L(m_j)$, $j=1,\ldots,r$.  The complex
% structure of $L(m_j)$ is such that $g\in S^1$ acts on $L(m_j)$ via
% complex multiplication with $m_j$.  The numbers $m_j\neq 0$,
% $j=1,\ldots,r$, are called rotation numbers.
% By the real splitting principle, they are defined
% also for even dimensional real oriented vector bundles, but in that
% case the $m_j$'s are well defined only up to an even number of sign
% changes.  Choose $(m_j)_{j=1,\ldots,r}$ a system of rotation
% numbers for $V/V^{S^1}$  For $j=1,\ldots,r$ define
% the numbers $r_j$ and $q_j$ by
%    $$m_j=n q_j + r_j \period$$
% For $0\leq k \leq \frac{n}{2}$ define
%   $$ I_k = \{j \in 1,\ldots,r\ |\ r_j= k \mbox{ or } n-k \}\period$$
% Notice that for $0<k\leq \frac{n}2$, $I_k$ contains exactly the
% indices of the rotation numbers of $V(k)$.
% For $k=0$, $I_0$ contains the indices of the rotation
% numbers corresponding to $V^{\Z_n}/V^{S^1}$.
% Define
%   $$ I_K = \cup_{0<k<\frac{n}2} I_k  \period$$

\vs
\bc
{\sc 5.2. Rigidity}
\ec
\vs

As in the beginning of Section~\ref{section-pushforwards},
let $\cE=\C/\Ll$
be an elliptic curve over $\C$ together with a 
2-torsion point and a generator of the cotangent space
to $\cE$ at zero.  We saw that we can canonically associate
to this data a double cover $\tilde{\cE}$ of $\cE$, and
the Jacobi sine function $s:\cE\ra\C$.

Let $X$ be a compact spin
$S^1$-manifold, i.e.\ a spin manifold with an $S^1$ action which
preserves the spin structure. Then the map $\pi:X\ra *$ is spin
oriented, hence it satisfies the hypothesis of
Definition~\ref{pushforward}.  Therefore we get a Grojnowski
pushforward $\pi^E_!:\E^*_{S^1}(X)^{[\pi]} \ra \E^*_{S^1}(*)=\cO_\cE$.

We will see that the rigidity phenomenon amounts to finding a global
(Thom) section in the sheaf $\E^*_{S^1}(X)^{[\pi]}$.  Since $s(x)$
is not a well-defined function on $\cE$, we cannot expect to find
such a global section on $\cE$.  However, if we take the pullback of the
sheaf $\E^*_{S^1}(X)^{[\pi]}$ along the covering map
$\tilde{\cE} \ra \cE$, we can show that the new sheaf has a global
section.

\vs
\noindent
\dfn{Convention.} From this point on, all the sheaves $\cF$ involved
will be considered over $\tilde{\cE}$, i.e.\ we will replace them
by the pullback of $\cF$ via the map $\tilde{\cE} \ra \cE$.
\vs

For our purposes, however, we need a more general version of
$\E^*_{S^1}(X)^{[\pi]}$, which involves a vector bundle $V\ra X$.
Consider now $V\ra X$ a spin $S^1$-vector bundle over a finite
$S^1$-$CW$ complex.

\bdefn \label{V-twisted-glue}
As in Definition~\ref{twisted-glue}, we define
$\phi_{\aa\bb}^{[V]}$ as the composition of three maps, where
the second one is multiplication by 
$\ll_{\aa\bb}^{[V]}=e^E_{S^1}(V^\aa/V^\bb)^{-1}$.  The bundle
$V^\aa/V^\bb = V^{\Z_n}/V^{S^1}$ is oriented as in the
previous subsection.

We then denote by $\E^*_{S^1}(X)^{[V]}$ the sheaf obtained
by gluing the sheaves $\cF_\aa$ defined in~\ref{local-sheaf},
using the twisted gluing functions $\phi_{\aa\bb}^{[V]}$.
\edefn

Notice that, if we take the map $f:X\ra *$ and $V=TX$,
we have $\E^*_{S^1}(X)^{[V]}=\E^*_{S^1}(X)^{[f]}$.
We now proceed to proving THEOREM C.

\bthm \label{thom-section}
If $\; V\ra X$ is a spin $S^1$-vector bundle over a finite
$S^1$-$CW$ complex, then the element 1 in the stalk of
$\E^*_{S^1}(X)^{[V]}$ at zero extends to a global section
on $\tilde{\cE}$, called the \dfn{Thom section}.
\ethm

\bp
To simplify notation, we are going to identify $\tilde{\cE}$ with
$\C/\tilde{\Ll}$, where $\tilde{\Ll} = \Z\oo_1+2\Z\oo_2$ is
the ``doubled'' lattice described in Section~\ref{section-pushforwards}.
We want now to think of points in $\tilde{\cE}$ as
points in $\C$, and of $\E^*_{S^1}(X)$ as the pullback
of $\E^*_{S^1}(X)$
on $\C$ via $\C \ra \C/\tilde{\Ll}$.  Then we call $\aa\in\C$
a torsion point if there is an integer $n>0$ such that
$n\aa\in\Ll$ (notice that torsion points are defined in terms
of $\Ll$, and not $\tilde{\Ll}$).  The smallest such $n$ is called the
exact order of $\aa$.  From Proposition~\ref{Jacobi-sine} b), we
know that if $a\in\Ll$, $s(x+a)=\pm s(x)$.  Since $n\aa\in\Ll$,
define $\eps=\pm 1$ by
        $$s(x+n\aa)=\eps s(x) \period$$

Now $\E^*_{S^1}(X)^{[V]}$ was obtained by gluing the sheaves
$\cF_\aa$ along the adapted open cover $(U_\aa)_\aa$.
So to give a global section $\mu$ of $\E^*_{S^1}(X)^{[V]}$ is the
same as to give global sections $\mu_\aa$ of  $\cF_\aa$ such that they
glue, i.e.\ $\phi_{\aa\bb}^{[V]}\mu_\aa = \mu_\bb$ for any $\aa$
and $\bb$ with $U_\aa \cap U_\bb \neq \emptyset$.  
From Definition~\ref{V-twisted-glue}, to give $\mu$ is the same as 
to give $\mu_\aa \in \ho^*_{S^1}(X^\aa)$ so that 
$t^*_{\bb-\aa}(i^*\mu_\aa \cdot e^E_{S^1}(V^\aa/V^\bb)^{-1}) = \mu_\bb$,
or $i^*\mu_\aa \cdot e^E_{S^1}(V^\aa/V^\bb)^{-1} = t^*_{\bb-\aa}\mu_\bb$
($i$ the inclusion $X^\bb\inc X^\aa$).  Because $\mu$ is supposed
to globalize 1, we know that
$\mu_0=1$.  This implies that $\mu_\bb=t^*_\bb e^E_{S^1}(V/V^\bb)^{-1}$
for $\bb$ in a small neighborhood of $0\in\C$.

In fact, we can show that this formula for $\mu_\bb$ is valid
for all $\bb\in \C$, as long as $\bb$ is not special. 
This means we have to check that
$\mu_\bb=t^*_\bb e^E_{S^1}(V/V^\bb)^{-1}$ exists in
$\ho^*_{S^1}(X^\bb)$ as long as $\bb$ is not
special.  $\bb$ not special means $X^\bb=X^{S^1}$.  Then consider
the bundle $V/V^{S^1}$.  We saw in the previous subsection
that according to the splitting principle, when pulled back
on the flag manifold, $V/V^{S^1}$ decomposes into a direct sum 
of line bundles
$L(m_1)\oplus \cdots \oplus L(m_r)$, where $m_j$ are the
rotation numbers.  The complex structure on $L(m)$ is such
$g\in S^1$ acts on $L(m)$ by complex multiplication with $g^m$.

Since $X^{S^1}$ is fixed by the $S^1$ action, we can apply
Proposition~\ref{formula-Chern-roots} in the Appendix:  Let
$x_j$ be the equivariant Chern root of $L(m_j)$, and $w_j$
its usual (nonequivariant) Chern root.  Then
$x_j=w_j+m_j u$, with $u$ the generator of $H^*(BS^1)$.
Therefore
$t^*_\bb e^E_{S^1}(V/V^\bb) = \prod_j t^*_\bb s(x_j) =
\prod_j t^*_\bb s(w_j+m_j u) = \prod_j s(w_j+m_j u + m_j\bb) =
\prod_j s(x_j + m_j\bb)$.

So we have
  $$ \mu_\bb = t^*_\bb e^E_{S^1}(V/V^{S^1})^{-1} =
           \prod_{j=1}^r s(x_j+m_j\bb)^{-1}  \period   $$
We show that $\mu_\bb$ belongs to $\ho^*_{S^1}(X^\bb)$
as long as $s(m_j\bb)\neq 0$ for all $j=1,\ldots,r$:
Since $V/V^{S^1}$ has only nonzero rotation numbers, it
has a complex structure.  But changing the orientations
of a vector bundle only changes the sign of the corresponding
Euler class, so in the formula above we can assume that
$V/V^{S^1}$ has a complex structure, for example the one
for which all $m_j>0$.  We group the $m_j$ which
are equal, i.e.\ for each $m>0$ we define the set of indices
$J_m = \{j \; | \; m_j=m \}$.  Now we get a
decomposition\footnote{This decomposition takes place
on $X^{S^1}$, while the decomposition into line bundles $L(m_j)$
takes place only on the flag manifold.} 
$V/V^{S^1}=\sum_{m>0}W(m)$, where $W(m)$ is the complex
$S^1$-vector bundle on which $g\in S^1$ acts by
multiplication with $g^m$.  Now we have to show that
$\prod_{j\in J_m}s(x_j+m\bb)^{-1}$ gives an element of
$\ho^*_{S^1}(X^\bb)$.  This would follow from 
Proposition~\ref{holomorphic} applied to the power series 
$Q(x)=s(x+m\bb)^{-1}$ and the vector bundle $W(m)$,
provided that $Q(x)$ is convergent.  But $s(x+m\bb)^{-1}$
is indeed convergent, since $s$ is meromorphic
on $\C$ and does not have a zero at $m\bb$.

Now we show that if $\bb$ is nonspecial,
 $s(m_j\bb)\neq 0$ for all $j=1,\ldots,r$:
Suppose $s(m_j\bb) = 0$.  Then $m_j\bb \in \Ll$, so $\bb$ is a 
torsion point, say of exact order $n$.  It follows that $n$ divides
$m_j$, which implies $X^{\Z_n}\neq X^{S^1}$.  But 
$X^\bb=X^{\Z_n}$, since $\bb$ has exact order $n$, so $X^\bb\neq X^{S^1}$
i.e.\ $\bb$ is special, contradiction.

So we only need to analyze what happens at a special point
$\aa\in\C$, say of exact order $n$.  We have to find a class 
$\mu_\aa \in \ho^*_{S^1}(X^\aa)$ such that
$\phi_{\aa\bb}^{[V]}\mu_\aa = \mu_\bb$, i.e.\ 
$t^*_{\bb-\aa}(i^*\mu_\aa \cdot e^E_{S^1}(V^\aa/V^\bb)^{-1})=
t^*_\bb e^E_{S^1}(V/V^\bb)^{-1}$.  Equivalently, we want a class
$\mu_\aa$ such that $i^* \mu_\aa = t^*_\aa e^E_{S^1}(V/V^\bb)^{-1} \cdot
e^E_{S^1}(V^\aa/V^\bb)$,
i.e.\ we want to lift the class $t^*_\aa e^E_{S^1}(V/V^\bb)^{-1} \cdot
e^E_{S^1}(V^\aa/V^\bb)$ from $\ho^*_{S^1}(X^\bb)$
to $\ho^*_{S^1}(X^\aa)$.  If we can do that, we are done, because the
class $(\mu_\aa)_{\aa\in\C}$ is a global section in
$\E^*_{S^1}(X)^{[V]}$, and it extends $\mu_0=1$ in the stalk at zero.
So it only remains to prove the following lemma, which is a 
generalization of the transfer formula of Bott and Taubes.
\ep

\blemma \label{transfer-formula}
Let $\aa$ be a special point of exact order $n$, and $V\ra X$ a spin
$S^1$-vector bundle.  Let $i: X^{S^1} \ra X^{\Z_n}$ be the inclusion
map.  Then there exists a class $\mu_\aa \in \ho^*_{S^1}(X^{\Z_n})$
such that
    $$i^* \mu_\aa = t^*_\aa e^E_{S^1}(V/V^{S^1})^{-1} \cdot
                e^E_{S^1}(V^{\Z_n}/V^{S^1}) \period$$
\elemma

\bp
We first study the class
$t^*_\aa e^E_{S^1}(V/V^{S^1})^{-1} \cdot e^E_{S^1}(V^{\Z_n}/V^{S^1})$
on each connected component of $X^{S^1}$ in $X^{\Z_n}$.  We
will see that it lifts naturally to a class on $X^{\Z_n}$.
The problem arises from the fact that we can have
two connected components of $X^{S^1}$ inside one connected
component of $X^{\Z_n}$, and in that case the two lifts will
differ by a sign.  We then show that the sign vanishes if
$V$ has a spin structure.

As in the previous subsection, let $N$ be a connected
component of $X^{S^1}$, and $P$ a connected component of
$X^{\Z_n}$ which contains $N$.
 
% Now we want to show that there exists a class 
% $\nu_P \in \ho^*_{S^1}(P)$  such that
% \begin{equation}
% i^* \nu_P = t^*_\aa e^E_{S^1}(V/V^{S^1})^{-1}_{\mid N} \cdot
%                 e^E_{S^1}(V^{\Z_n}/V^{S^1})_{\mid N} \period
% \end{equation}
We now calculate $t^*_\aa e^E_{S^1}(V/V^{S^1})^{-1}$, regarded
as a class on $N$.
From the decomposition~(\ref{VS})
$V/V^{S^1} = V^{\Z_n}/V^{S^1} \oplus V(K)_{\mid N}
\oplus V(\frac{n}{2})_{\mid N}$ and from the table,
we get the following formula:
\begin{eqnarray} \label{formula}
    t^*_\aa e^E_{S^1}(V/V^{S^1})^{-1} & = & 
         (-1)^\ss \cdot e^E_{S^1}(V/V^{S^1})^{-1}_{cx}  \nonumber \\
      &  = & (-1)^\ss \cdot \prod_{j\in I_0} s(x_j+m_j^* \aa)^{-1} \cdot 
         \prod_{j \in I_K} s(x_j+m_j^* \aa)^{-1} \cdot 
           \prod_{j\in I_{n/2}} s(x_j+m_j^* \aa)^{-1} 
\end{eqnarray}

Before we analyze each term in the above formula, recall that we
defined the number $\eps=\pm 1$ by $s(x+n\aa)=\eps s(x)$.

a) $j\in I_0$: Here we chose the complex structure 
$(V^{\Z_n}/V^{S^1})_{cx}$ such that all $m_j^*>0$.  Then,
since $s(x_j+m_j^* \aa)=s(x_j+q_j^* n \aa)=\eps^{q_j^*}s(x_j)$,
we have:
$ \prod_{j\in I_0} s(x_j+m_j \aa)^{-1} =
  \eps^{\sum_{I_0}q_j^*} \cdot \prod_{I_0} s(x_j)^{-1} =
  \eps^{\sum_{I_0}q_j^*} \cdot e^E_{S^1}(V^{\Z_n}/V^{S^1})^{-1}_{cx} =
  \eps^{\sum_{I_0}q_j^*} \cdot (-1)^{\ss(0)} \cdot
                     e^E_{S^1}(V^{\Z_n}/V^{S^1})^{-1}_{or}$.
So we get eventually
\begin{equation} 
 \prod_{j\in I_0} s(x_j+m_j^* \aa)^{-1} = 
   \eps^{\sum_{I_0}q_j^*} \cdot (-1)^{\ss(0)} \cdot 
                       e^E_{S^1}(V^{\Z_n}/V^{S^1})^{-1}_{or} \period
\end{equation}

b) $j\in I_K$, i.e.\ $j\in I_k$ for some $0<k<\frac{n}{2}$.
The complex structure on $V(k)$ is such that $g=e^{2\pi i/n}\in\Z_n$
acts by complex multiplication with $g^k$.  Notice that in
the previous subsection we defined the complex structure
on $V/V^{S^1}$ to come from the decompostion~(\ref{VS}).
This implies that $m_j^*=nq_j^*+k$, and therefore
$s(x_j+m_j^*\aa)=s(x_j+q_j^*n\aa+k\aa)=\eps^{q_j^*}s(x_j+k\aa)$.

Consider $\mu_k$ the equivariant class on $P$ 
corresponding to the complex vector bundle $V(k)$ with its chosen
complex orientation, and the convergent power series
$Q(x)=s(x+k \aa)^{-1}$.  Then $i^*\mu_k=\prod_{I_k} s(x_j+k \aa)^{-1}$.
Define $\mu_K=\prod_{0<k<\frac n2}\mu_k$.
Using the above formula for $s(x_j+m_j^*\aa)$ with $j\in I_k$, we
obtain
\begin{equation} 
\prod_{j\in I_K} s(x_j+m_j^* \aa)^{-1} = 
   \eps^{\sum_{I_K}q_j^*} \cdot
             (-\eps)^{\ss(K)} \cdot
                          i^*\mu_K \period
\end{equation}

c) $j\in I_{n/2}$.  The complex structure on $i^*V(\frac{n}{2})$
is the one for which all $m_j^*>0$.  The rotation
numbers satisfy $m_j^*=q_j^* n + \frac{n}{2}$, hence
$s(x_j+m_j^* \aa)=\eps^{q_j^*}s(x_j+\frac{n}{2}\aa)$.
Consider the power series $Q(x)=s(x+ \frac{n}{2}\aa)^{-1}$.
$Q(x)$ satisfies 
$Q(-x)=s(-x+ \frac{n}{2}\aa)^{-1}=-s(x- \frac{n}{2}\aa)^{-1}=
-\eps s(x+ \frac{n}{2}\aa)^{-1}=(-\eps)Q(x)$, so $Q(x)$ is either even
or odd.  According to Definition~\ref{defn-real-equiv-char-class},
since $V(\frac n2)_{or}$ is a real oriented even dimensional vector
bundle, $Q(x)$ defines a class $\mu_{\frac{n}{2}}=\mu_Q(V(\frac n2))$,
which is a clas on $P$.  Now from the table, $i^*V(\frac n2)_{or}$ and
$(i^*V(\frac{n}{2}))_{cx}$ differ by the sign $(-1)^{\ss(\frac n2)}$,
so Lemma~\ref{real-class} (with $\gg=-\eps$) implies that
$i^*\mu_{\frac{n}{2}} = (-\eps)^{\ss(\frac{n}{2})} 
             \prod_{j\in I_k} s(x_j+\frac{n}{2} \aa)^{-1}$.
Finally we obtain
\begin{equation} \label{half}
\prod_{j\in I_{n/2}} s(x_j+m_j^* \aa)^{-1} =
   \eps^{\sum_{I_{n/2}}q_j^*} \cdot 
           (-\eps)^{\ss(\frac{n}{2})} \cdot
                 i^*\mu_{\frac{n}{2}} \period
\end{equation}

Now, if we put together equations (\ref{formula})--(\ref{half})
and (\ref{signs}), and define $\mu_P:= \mu_K \cdot \mu_{\frac{n}{2}}$,
we have just proved that 
$t^*_\aa e^E_{S^1}(V/V^{S^1})^{-1} = \eps^{\ss(N)}
           \cdot e^E_{S^1}(V^{\Z_n}/V^{S^1})^{-1} \cdot i^*\mu_P$, or
\begin{equation} \label{lift-NP}
  t^*_\aa e^E_{S^1}(V/V^{S^1})^{-1} \cdot e^E_{S^1}(V^{\Z_n}/V^{S^1}) 
         = \eps^{\ss(N)} \cdot i^*\mu_P \comma
\end{equation}
where
  $$ \ss(N)=\sum_{I_0}q_j^* + \sum_{I_K}q_j^* +
          \sum_{I_{n/2}}q_j^* + \ss(K)+\ss(\frac{n}{2}) \period$$

Now we want to describe $\ss(N)$ in terms of the correct rotation
numbers $m_j$ of $V/V^{S^1}$.  Recall that $m_j$ are the same as
$m_j^*$ up to sign and a permutation.  Denote by $\equiv$ equality modulo
2.  We have the following cases:
\bi
\item[a)] $j\in I_0$.  Suppose $m_j=-m_j^*$.  Then $q_j=-q_j^*$, which
 implies $q_j^*\equiv q_j$.  Therefore $\sum_{I_0}q_j^* \equiv \sum_{I_0}q_j$. 
\item[b)] $j\in I_K$.  Let $0<k<\frac{n}{2}$.  Suppose
 $m_j=-m_j^* = -q_j^* n-k = -(q_j^*+1)n+(n-k)$.  Then $q_j=-q_j^*-1$, which
 implies $q_j^*+1\equiv q_j$.  So modulo 2, the sum $\sum_{I_K}q_j^*$
 differs  from $\sum_{I_K}q_j$ by the number of the sign differences
 $m_j=-m_j^*$.  But by definition of rotation numbers, the number of sign
 differences in two systems of rotation numbers is precisely the sign
 difference $\ss(K)$ between the two corresponding orientations
 of $i^*V(K)$.  Therefore, $\sum_{I_K}q_j^* + \ss(K) \equiv \sum_{I_K}q_j$.
\item[c)] $j\in I_{n/2}$.  Suppose 
 $m_j=-m_j^*=-q_j^* n-\frac n2=-(q_j^*+1)n+\frac n2$.  Then this implies
 $q_j^*+1\equiv q_j$, so by the same reasoning as in b)
 $\sum_{I_{n/2}}q_j^* + \ss(\frac n2) \equiv \sum_{I_{n/2}}q_j$.
\ei
We finally get the following formula for $\ss(N)$
  $$\ss(N)\equiv\sum_{I_0}q_j + \sum_{I_K}q_j 
    +\sum_{I_{n/2}}q_j \period$$

In the next lemma we will show that, for $N$ and $\tilde{N}$ two different
connected components of $X^{S^1}$ inside $P$, $\ss(N)$ and
$\ss(\tilde{N})$ are congruent modulo 2, so the class
$\eps^{\ss(N)} \cdot \mu_P$ is well-defined, i.e.\ independent of $N$.
Now recall that $P$ is a connected component of $X^{\Z_n}$.  Therefore
$\ho^*_{S^1}(X^{\Z_n})=\oplus_P\ho^*_{S^1}(P)$, so we can define
  $$\mu_\aa:=\sum_P \eps^{\ss(N)}\cdot\mu_P \period$$
This is a well-defined class in $\ho^*_{S^1}(X^{\Z_n})$, so by equation
(\ref{lift-NP}), Lemma~\ref{transfer-formula} is finally proved.
\ep

\blemma
In the conditions of the previous lemma, $\ss(N)$ and 
$\ss(\tilde{N})$ are congruent modulo 2.
\elemma

\bp
The proof follows Bott and Taubes~\cite{BT}.
Denote by $S^2(n)$ the $2$-sphere with the $S^1$-action which 
rotates $S^2$ $n$ times around the north-south axis as we go once 
around $S^1$.  Denote by $N^+$ and $N^-$ its North and South poles,
respectively.  Consider a path in $P$ which connects $N$ with $\tilde{N}$,
and touches $N$ or $\tilde{N}$ only at its endpoints.  By rotating
this path with the $S^1$-action, we obtain a subspace of $P$ which
is close to being an embedded $S^2(n)$.  Even if it is not, we 
can still map equivariantly $S^2(n)$ onto this rotated path.  
Now we can pull back the bundles from $P$ to $S^2(n)$ (with their
correct orientations).  The rotation numbers are the same, since
the North and the South poles are fixed by the $S^1$-action, as are
the endpoints of the path.  

Therefore we have translated the problem to the case when we have
the $2$-sphere $S^2(n)$ and corresponding bundles over it, and we
are trying to prove that $\ss(N^+)\equiv\ss(N^-)$ modulo 2.
The only problem
would be that we are not using the whole of $V$, but only 
$V/V^{S^1}$.  However, the difference between these two bundles is
$V^{S^1}$, whose rotation numbers are all zero, so they do not
influence the result.

Now Lemma 9.2 of~\cite{BT} says that any even-dimensional oriented
real vector bundle $W$ over $S^2(n)$ has a complex structure.  In
particular, the pullbacks of $V^{S^1}$, $V(K)$, and $V(\frac n2)$
have complex structure, and the rotation numbers can be chosen to
be the $m_j$ described above.  Say the rotation numbers at the 
South pole are $\tilde{m}_j$ with the obvious notation conventions.
Then Lemma 9.1 of~\cite{BT} says that, up to a permutation, 
$m_j -\tilde{m}_j=n(q_j-\tilde{q}_j)$, and 
$\sum q_j \equiv \sum \tilde{q}_j$ modulo 2.  But this means that
$\ss(N^+)\equiv\ss(N^-)$ modulo 2,
i.e.\ $\ss(N)\equiv \ss(\tilde{N})$ modulo 2.
\ep

\bcor (The Rigidity theorem of Witten) \label{rigidity-theorem}
If $X$ is a spin manifold with an $S^1$-action, then the equivariant 
elliptic genus of $X$ is rigid i.e.\ it is a constant power series. 
\ecor

\bp
By lifting the $S^1$-action to a double cover of 
$S^1$, we can make the $S^1$-action preserve the spin structure.
Then with this action $X$ is a spin $S^1$-manifold.

At the beginning of this Section, we say that if $X$ is a compact 
spin $S^1$-manifold, i.e.\ the map $\pi:X \ra *$ is 
spin, then we have the Grojnowski pushforward, which 
is a map of sheaves 
  $$ \pi^E_!: \E^*_{S^1}(X)^{[\pi]} \ra 
                            \E^*_{S^1}(*)=\cO_\cE \period$$
The Grojnowski pushforward $\pi^E_!$, if we 
consider it at the level of stalks at $0 \in \cE$, is nothing
but the elliptic pushforward in $\ho^*_{S^1}$-theory, as described in 
Corollary \ref{local-pushforward}.  So consider the element 
$1$ in the stalk at $0$ of the sheaf
$\E^*_{S^1}(X)^{[\pi]}=\E^*_{S^1}(X)^{[TX]}$.  

From Theorem~\ref{thom-section}, since $TX$ is spin, 1 extends
to a global section of $\E^*_{S^1}(X)^{[TX]}$.  Denote 
this global section by boldface $\one$.  Because $\pi^E_!$ is a 
map of sheaves, it follows that $\pi^E_!(\one)$ is a global 
section of $\E^*_{S^1}(*) = \cO_\cE$,
i.e.\ a global holomorphic function on the elliptic curve $\cE$.
But any such function has to be constant.  This means that
$\pi^E_!(1)$, which is the equivariant elliptic genus of $X$,
extends to $\pi^E_!(\one)$, which is constant.  This is precisely
equivalent to the elliptic genus being rigid.
\ep

The extra generality we had in Theorem~\ref{thom-section} allows
us now to extend the Rigidity theorem to families of elliptic
genera.  This was stated as THEOREM D in Section~\ref{section-results}.

\bthm \label{rigidity-families}
(Rigidity for families) Let $F \ra E\llra{\pi} B$ be an 
$S^1$-equivariant fibration such that the fibers are spin 
in a compatible way, i.e.\ the projection map $\pi$ is spin oriented.
Then the elliptic genus of the family, which is 
$\pi^E_!(1)\in H^{**}_{S^1}(B)$, is constant as a rational
function in $u$, i.e.\ if we invert $u$.
\ethm

\bp
We know that the map
   $$ \pi^E_!:\E^*_{S^1}(E)^{[\pi]} \ra \E^*_{S^1}(B) $$
when regarded at the level of stalks at zero is the usual
equivariant elliptic pushforward in $\ho^*_{S^1}(\blank)$.
Now $\pi^E_!(1)\in \ho^*_{S^1}(B)$ is the elliptic genus of 
the family.  We have
$\E^*_{S^1}(E)^{[\pi]}\cong \E^*_{S^1}(E)^{[\tau(F)]}$,
where $\tau(F) \ra E$ is the bundle of tangents along the fiber.

Since $\tau(F)$ is spin, Theorem~\ref{thom-section} allows us to extend
1 to the Thom section $\one$.  Since $\pi^E_!$ is a map of
sheaves, it follows that $\pi^E_!(1)$, which is the elliptic genus 
of the family, extends to a global section in $\E^*_{S^1}(B)$.  
So, if $i:B^{S^1}\inc B$ is the inclusion of 
the fixed point submanifold in B, $i^* \pi^E_!(\one)$ gives a 
global section in $\E^*_{S^1}(B^{S^1})$.
Now this latter sheaf is free as a sheaf of $\cO_\cE$-modules,
so any global section is constant. But 
$i^*:\ho^*_{S^1}(B)\ra \ho^*_{S^1}(B^{S^1})$ is an isomorphism
if we invert $u$.
\ep

We saw in the previous section that, if $f: X\ra Y$ is an $S^1$-map
of compact $S^1$-manifolds such that the restrictions 
$f:X^\aa \ra Y^\aa$ are oriented maps, we have the Grojnowski pushforward
    $$f^E_!:\E^*_{S^1}(X)^{[f]} \ra \E^*_{S^1}(Y) \period$$
Also, in some cases, for example when $f$ is a spin $S^1$-fibration,
we saw that $\E^*_{S^1}(X)^{[f]}$ admits a Thom section.  This raises
the question if we can describe $\E^*_{S^1}(X)^{[f]}$ as $\E^*_{S^1}$
of a Thom space.  It turns out that, up to a line bundle over $\cE$
(which is itself $\E^*_{S^1}$ of a Thom space), this indeed happens:

Let $f: X\ra Y$ be an $S^1$-map as above.  Embed $X$ into an 
$S^1$-representation $W$, $i:X\inc W$.  ($W$ can be also thought as an
$S^1$-vector bundle over a point.)  Look at the embedding 
$f\times i: X \inc Y\times W$.  Denote by $V=\nu(f)$, the normal
bundle of $X$ in this embedding (if we were not in the equivariant
setup, $\nu(f)$ would be the stable normal bundle to the map $f$).

\bprop \label{identification-thom-sheaf}
With the previous notations,
  $$\E^*_{S^1}(X)^{[f]}\cong 
           \E^*_{S^1}(DV,SV)\ts\E^*_{S^1}(DW,SW)^{-1}\comma$$
where $DV$, $SV$ are the disk and the sphere bundles of $V$, respectively.
\eprop

\bp
From the embedding $X\inc Y\times W$, we have the following isomorphism
of vector bundles:
  $$TX\oplus V\cong f^*TY\oplus W\period $$
So, in terms of $S^1$-equivariant elliptic Euler classes we have
$e^E_{S^1}(V^\aa/V^\bb)=e^E_{S^1}(X^\aa/X^\bb)^{-1}\cdot
     f^*e^E_{S^1}(Y^\aa/Y^\bb)\cdot e^E_{S^1}(W^\aa/W^\bb)$.
Rewrite this as 
   $$\ll_{\aa\bb}^{[f]} =
        e^E_{S^1}(V^\aa/V^\bb)\cdot e^E_{S^1}(W^\aa/W^\bb)^{-1} \comma$$
where $\ll_{\aa\bb}^{[f]}$ is the twisted cocycle from
Definition~\ref{twisted-glue}.

Notice that we can extend Definition~\ref{V-twisted-glue} to virtual
bundles as well.  In other words, we can define
$\E^*_{S^1}(X)^{[-V]}$ to be $\E^*_{S^1}(X)$ twisted by the cocycle
$\ll_{\aa\bb}^{[-V]}=e^E_{S^1}(V^\aa/V^\bb)$.  The above formula
then becomes
  $$\ll_{\aa\bb}^{[f]} = \ll_{\aa\bb}^{[-V]} \cdot
                      \ll_{\aa\bb}^{[W]}           \comma$$
which implies that
\begin{equation} \label{formula-twisted-sheaf}
     \E^*_{S^1}(X)^{[f]} = 
            \E^*_{S^1}(X)^{[-V]}\ts \E^*_{S^1}(X)^{[W]} \period
\end{equation}
So the proposition is finished if we can show that
for a general vector bundle $V$
   $$\E^*_{S^1}(DV,SV)=\E^*_{S^1}(X)^{[-V]} \period $$
Indeed, multiplication by the equivariant elliptic Thom classes
on each stalk gives the following commutative diagram, where
the rows are isomorphisms:
 $$\xymatrix{ 
    H^*_{S^1}(X^\aa)\ts_{\C[u]}\cO_\cE(U-\aa) 
            \ar[d]_{e^E_{S^1}(V^\aa/V^\bb)\cdot i^*} 
                \ar[rr]^{\cdot t^*_\aa \phi^E_{S^1}(V^\aa)} & &
    H^*_{S^1}(DV^\aa,SV^\aa)\ts_{\C[u]}\cO_\cE(U-\aa)
            \ar[d]^{i^*}                                          \\  
    H^*_{S^1}(X^\bb)\ts_{\C[u]}\cO_\cE(U-\aa)
            \ar[d]_{t^*_{\bb-\aa}}
                \ar[rr]^{\cdot t^*_\aa \phi^E_{S^1}(V^\bb)} & &
    H^*_{S^1}(DV^\bb,SV^\bb)\ts_{\C[u]}\cO_\cE(U-\aa)
            \ar[d]^{t^*_{\bb-\aa}}                                 \\
    H^*_{S^1}(X^\bb)\ts_{\C[u]}\cO_\cE(U-\bb) 
            \ar[rr]^{\cdot t^*_\bb \phi^E_{S^1}(V^\bb)} & &
    H^*_{S^1}(DV^\bb,SV^\bb)\ts_{\C[u]}\cO_\cE(U-\bb)     \period
            }    $$
Notice that $\E^*_{S^1}(DW,SW)$ is an invertible sheaf, because it is
the same as the structure sheaf $\E^*_{S^1}(*)=\cO_{\cE}$ twisted
by the cocycle $\ll_{\aa\bb}^{[W]}$.  In fact, we can identify it by
the same method we used in Proposition~\ref{example-thom-sheaf}.
\ep

In the language of equivariant spectra (see Chapter 8 of~\cite{IJ})
we can say more:  With the notation we used in 
Proposition~\ref{identification-thom-sheaf},
we define a virtual vector bundle $Tf$, the tangents along the fiber,
by 
   $$TX = Tf \oplus f^*TY\period$$
Using the formula $TX\oplus V= f^*TY\oplus W$, it follows
that $-Tf = V \ominus W$.  From equation~(\ref{formula-twisted-sheaf})
it follows that
    $$\E^*_{S^1}(X)^{[f]} = \tilde{\E}^*_{S^1}(X^{-Tf}) \comma$$
where $\tilde{\E}^*_{S^1}$ is reduced cohomology, and
$X^{-Tf}$ is the $S^1$-equivariant spectrum obtained by
the Thom space of $V$ desuspended by $W$.

\vs

\appendix

\section{Equivariant characteristic classes} 
                               \label{section-char-classes}

The results of this section are well-known, with the exception
of the holomorphicity result Proposition~\ref{holomorphic}.

Let $V$ be a complex $n$-dimensional $S^1$-equivariant vector 
bundle over an $S^1$-$CW$ complex $X$.  Then to any power series 
$Q(x)\in\C[\![x]\!]$ starting with 1 we are going to associate 
by Hirzebruch's formalism (see~\cite{Hi}) a multiplicative 
characteristic class $\mu_Q(V)_{S^1} \in H^{**}_{S^1}(X)$.  
(Recall that $H^{**}_{S^1}(X)$ is the completion of $H^*_{S^1}(X)$.)

Consider the Borel construction for both $V$ and $X$:
$V_{S^1} = V \times_{S^1} ES^1 \ra  X \times_{S^1} ES^1=X_{S^1}$.
$V_{S^1}\ra X_{S^1}$ is a complex vector bundle over a paracompact
space, hence we have a classifying map $f_V: X_{S^1} \ra BU(n)$.
We define $c_j(V)_{S^1}$, the equivariant $j$'th Chern class of $V$,
as the image via $f_V^*$ of the universal $j$'th Chern
class $c_j\in H^*BU(n)=\C[c_1,\ldots,c_n]$.  Now look at the product
$Q(x_1)Q(x_2)\cdots Q(x_n)$.  It is a power series in 
$x_1,\ldots,x_n$ which is symmetric under permutations of the
$x_j$'s, hence it can be expressed as another power series in 
the elementary symmetric functions 
$\ss_j=\ss_j(x_1,\ldots,x_n)$:
   $$Q(x_1)\cdots Q(x_n) = P_Q(\ss_1,\ldots,\ss_n) \period$$
Notice that $P_Q(c_1,\ldots,c_n)$ lies not in 
$H^*BU(n)$, but in its completion $H^{**}BU(n)$.
The map $f_V^*$ extends to a map $H^{**}BU(n)\ra H^{**}(X_{S^1})$.

\bdefn \label{defn-cx-equiv-char-class}
Given the power series $Q(x)\in \C[\![x]\!]$ and the 
complex $S^1$-vector bundle $V$ over $X$, there is a canonical
complex equivariant characteristic class 
$\mu_Q(V)_{S^1}\in H^{**}(X_{S^1})$,
given by 
    $$\mu_Q(V)_{S^1} := P_Q(c_1(V)_{S^1},\ldots,c_n(V)_{S^1})
      =f^*_V P_Q(c_1,\ldots,c_n) \period$$
\edefn

\bremark \label{complex-chern-roots}
If $T^n\inc BU(n)$ is a maximal torus, then then 
$H^*BT^n=\C[x_1,\ldots,x_n]$, and the $x_j$'s are called 
the universal Chern roots.  The map 
$H^*BU(n) \ra H^*BT^n$ is injective, and its image can
be identified as the Weyl group invariants of $H^*BT^n$.
The Weyl group of $U(n)$ is the symmetric group on $n$
letters, so $H^*BU(n)$ can be identified 
as the subring of symmetric polynomials in $\C[x_1,\ldots,x_n]$.
Similarly, $H^{**}BU(n)$ is the subring of symmetric power series
in $\C[\![x_1,\ldots,x_n]\!]$.
Under this interpretation, $c_j=\ss_j(x_1,\ldots,x_n)$.
This allows us to identify $Q(x_1) \cdots Q(x_n)$ with the element
$P_Q(c_1,\ldots,c_n) \in H^{**}BU(n)$.
\eremark

\bdefn \label{equivariant-Chern-roots}
We can write formally $\mu_Q(V)_{S^1}=Q(x_1) \cdots Q(x_n)$.
$x_1, \ldots, x_n$ are called the equivariant Chern roots of $V$.
\edefn

Here is a standard result about the equivariant Chern roots:

\bprop \label{formula-Chern-roots}
Let $V(m)\ra X$ be a complex $S^1$-vector bundle such that
the action of $S^1$ on $X$ is trivial.  Suppose that $g\in S^1$
acts on $V(m)$ by complex multiplication with $g^{m}$.
If $x_i$ are the equivariant Chern roots of $V(m)$,
and $w_i$ are its usual (nonequivariant) Chern roots,
then
   $$ x_i = w_i + mu \comma$$
where $u$ is the generator of $H^*_{S^1}(*)=H^*BS^1$.
\eprop

We want now to show that the class we have just constructed,
$\mu_Q(V)_{S^1}$, is holomorphic in a certain sense, provided
$Q(x)$ is the expansion of a holomorphic function around zero.
But first, let us state a classical lemma in the theory of 
symmetric functions.

\blemma \label{symmetric-holo}
Suppose $Q(y_1,\ldots,y_n)$ is a holomorphic (i.e.\ convergent) power
series, which is symmetric under permutations of the $y_j$'s.
Then the power series $P_Q$ such that
   $$Q(y_1,\ldots,y_n)=P_Q(\ss_1(y_1,\ldots,y_n),\ldots,
             \ss_n(y_1,\ldots,y_n))\comma$$
is holomorphic.
\elemma

We have mentioned above that $\mu_Q(V)_{S^1}$ belongs to 
$H^{**}_{S^1}(X)$.  This ring 
is equivariant cohomology tensored with power series.  It 
contains $\ho^*_{S^1}(X)$ as a subring, corresponding to the
holomorphic power series.

\bprop \label{holomorphic}
If $Q(x)$ is a convergent power series, then $\mu_Q(V)_{S^1}$ 
is a holomorphic class, i.e.\ it belongs to the subring 
$\ho^*_{S^1}(X)$ of $H^{**}_{S^1}(X)$.
\eprop

\bp 
We have $\mu_Q(V)_{S^1} = P(c_1(V)_{S^1},\ldots,c_n(V)_{S^1})$,
where we write $P$ for $P_Q$.

Assume $X$ has a trivial $S^1$-action.  It is easy to see that
$H^*_{S^1}(X)=(H^0(X)\ts_\C \C[u]) \oplus \textrm{nilpotents}$.
Hence we can write $c_j(E)_{S^1}=f_j+\aa_j$, with
$f_j\in H^0(X)\ts_\C \C[u]$, and $\aa_j$ nilpotent in
$H^*_{S^1}(X)$.  We expand $\mu_Q(V)_{S^1}$ in Taylor expansion
in multiindex notation.  We make the following notations: 
$\ll =(\ll_1,\cdots,\ll_n) \in\N^n$, 
$|\ll|=\ll_1 + \cdots + \ll_n$, and 
$\aa^\ll=\aa_1^{\ll_1}\cdots\aa_n^{\ll_n}$.
Now we consider the Taylor expansion of $\mu_Q(V)_{S^1}$ 
in multiindex notation:
  $$\mu_Q(V)_{S^1} = P(\ldots,c_j(V)_{S^1},\ldots)=
     \sum_\ll \frac{\del^{|\ll|} P}{\del c^\ll}
          (\ldots,f_j,\ldots) \cdot \aa^\ll \period$$
This is a finite sum, since $\aa_j$'s are nilpotent.
We want to show that $\mu_Q(V)_{S^1}\in \ho^*_{S^1}(X)$.
$\aa^\ll$ lies in $\ho^*_{S^1}(X)$, since it lies even
in $H^*_{S^1}(X)$.  So we only have to show that 
$\frac{\del^{|\ll|} P}{\del c^\ll}(\ldots,f_j,\ldots)$
lies in $\ho^*_{S^1}(X)$.

But $f_j\in H^0(X)\ts_\C \C[u] = \C[u]\oplus\cdots\oplus\C[u]$,
with one $\C[u]$ for each connected component of $X$.  If we 
fix one such component $N$, then the corresponding component 
$f_j^{(N)}$ lies in $\C[u]$.  According to Lemma~\ref{symmetric-holo},
$P$ is holomorphic around $(0,\ldots,0)$, hence so is
$\frac{\del^{|\ll|} P}{\del c^\ll}$.  Therefore 
$\frac{\del^{|\ll|} P}{\del c^\ll}
(\ldots,f_j^{(N)}(u),\ldots)$ is holomorphic in $u$ around $0$,
i.e.\ it lies in $\cO_{\C,0}$.  Collecting the terms for the 
different connected components of $X$, we finally get 
  $$\frac{\del^{|\ll|} P}{\del c^\ll}(\ldots,f_j,\ldots)
         \in \cO_{\C,0}\oplus\cdots\oplus\cO_{\C,0}=
             H^0(X)\ts_\C \cO_{\C,0} \period$$
But $H^0(X)\ts_\C \cO_{\C,0}\subseteq H^*(X)\ts_\C \cO_{\C,0}
   = H^*_{S^1}(X) \ts_{\C[u]}\cO_{\C,0}=\ho^*_{S^1}(X)$,
so we are done.

If the $S^1$-action on $X$ is not trivial, look at the following
exact sequence associated to the pair $(X,X^{S^1})$:
$$ 0 \ra T \inc H^{*}_{S^1}(X) \llra{i^{*}} H^{*}_{S^1}(X^{S^1}) 
        \llra{\delta} H^{*+1}_{S^1}(X,X^{S^1}) \comma$$
where $T$ is the torsion submodule of $H^*_{S^1}(X)$.  (The fact
that $T=\ker i^*$ follows from the following arguments:
on the one hand, $\ker i^*$ is torsion, because of the localization
theorem; on the other hand, $H^*_{S^1}(X^{S^1})$ is free, hence
all torsion in $H^*_{S^1}(X)$ maps to zero via $i^*$.)
Also, since $T$ is a direct sum of torsion modules of the form
$\C[u]/(u^n)$
  $$T \ts_{\C[u]}\cO_{\C,0}\cong T \cong T \ts_{\C[u]}\C[\![u]\!]\period$$

Now tensor the above exact sequence with $\cO_{\C,0}$ and $\C[\![u]\!]$
over $\C[u]$:
  $$\xymatrix{  
       0 \ar[r] & 
         T\; \ar@{^{(}->}[r] \ar@{=}[d] &
            \ho^*_{S^1}(X) \ar[r]^{i^*} \ar@{^{(}->}[d]^s &
                \ho^*_{S^1}(X^{S^1}) \ar[r]^\delta \ar@{^{(}->}[d]^t &
                    \ho^{*+1}_{S^1}(X,X^{S^1}) \ar@{=}[d]  \\
       0 \ar[r] & 
         T\; \ar@{^{(}->}[r]&
            H^{**}_{S^1}(X) \ar[r]^{i^{*}} &
                H^{**}_{S^1}(X^{S^1}) \ar[r]^\delta &
                    H^{**+1}_{S^1}(X,X^{S^1}) \period        
               }$$
We know $\aa:= \mu_Q(V)_{S^1} \in H^{**}_{S^1}(X)$.  Then
$\bb:=i^*\mu_Q(V)_{S^1}=i^*\aa$ was shown previously to be in the 
image of $t$, i.e.\ $\bb = t\tilde{\bb}$.  
$\delta\bb=\delta i^*\aa=0$, so $\delta t \tilde{\bb}=0$, hence
$\delta \tilde{\bb}=0$.  Thus $\tilde{\bb}\in \im\; i^*$, so
there is an $\tilde{\aa}\in\ho^{*}_{S^1}(X)$ such that
$\tilde{\bb}=i^*\tilde{\aa}$.  $s\tilde{\aa}$ might not
equal $\aa$, but $i^*(\aa-\tilde{\aa})=0$, so 
$\aa-\tilde{\aa}\in T$.  Now,
$\tilde{\aa}+(\aa-\tilde{\aa}) \in \ho^{*}_{S^1}(X)$,
and $s(\tilde{\aa}+(\aa-\tilde{\aa})=\aa$, which 
shows that indeed $\aa\in \im\; s = \ho^{*}_{S^1}(X)$.
\ep

There is a similar story when $V$ is an oriented
$2n$-dimensional real $S^1$-vector bundle over a 
finite $S^1$-$CW$ complex $X$.  We classify $V_{S^1} \ra X_{S^1}$
by a map $f_V: X_{S^1} \ra BSO(2n)$.
$H^*BSO(2n)=\C[p_1,\ldots,p_n]/ (e^2-p_n)$, where $p_j$ and $e$
are the universal Pontrjagin and Euler classes, respectively.
The only problem now is that in order to define characteristic
classes over $BSO(2n)$ we need the initial power series 
$Q(x)\in \C[\![x]\!]$ to be either even or odd:

\bremark \label{real-chern-roots}
As in Remark~\ref{complex-chern-roots}, if $T^n\inc BSO(2n)$ 
is a maximal torus, then the map $H^*BSO(2n) \ra H^*BT^n$
is injective, and its image can be identified as the Weyl group
invariants of $H^*BT^n$.  Therefore $H^*BSO(2n)$ can be thought of
as the subring of symmetric polynomials in $\C[x_1,\ldots,x_n]$
which are invariant under an even number of sign changes of the $x_j$'s.
A similar statement holds for $H^{**}BSO(2n)$.
Under this interpretation, $p_j=\ss_j(x_1^2,\ldots,x_n^2)$ and
$e=x_1\cdots x_n$.

So, if we want $Q(x_1) \cdots Q(x_n)$ to be interpreted as an element of
$H^{**}BSO(2n)$, we need to make it invariant under an even number of 
sign changes.  But this is clearly true if $Q(x)$ is either
an even or an odd power series.

Let us be more precise:
\bi
\item[a)] $Q(x)$ is even, i.e.\ $Q(-x)=Q(x)$.  Then there is another
power series $S(x)$ such that $Q(x)=S(x^2)$, so 
$Q(x_1) \cdots Q(x_n)=S(x_1^2) \cdots S(x_n^2)=
P_S(\ldots,\ss_j(x_1^2,\ldots,x_n^2),\ldots)=
P_S(\ldots,p_j,\ldots)$.
\item[b)] $Q(x)$ is odd, i.e.\ $Q(-x)=-Q(x)$.  Then there is another
power series $R(x)$ such that $Q(x)=x T(x^2)$, so
$Q(x_1) \cdots Q(x_n)=x_1\cdots x_n T(x_1^2) \cdots T(x_n^2)=
x_1\cdots x_n P_T(\ldots,\ss_j(x_1^2,\ldots,x_n^2),\ldots)=
e\cdot P_T(\ldots,p_j,\ldots)$. 
\ei
\eremark

\bdefn \label{defn-real-equiv-char-class}
Given the power series $Q(x)\in\C[\![x]\!]$ which is either even or
odd, and the real oriented $S^1$-vector bundle $V$ over $X$, there
is a canonical real equivariant characteristic class 
$\mu_Q(V)_{S^1}\in H^{**}_{S^1}(X)$, defined by pulling back 
the element $Q(x_1) \cdots Q(x_n)\in H^{**}BSO(2n)$ via the 
classifying map $f_V: X_{S^1} \ra BSO(2n)$.
\edefn

Proposition~\ref{holomorphic} can be adapted to show that, 
if $Q(x)$ is a convergent power series, $\mu_Q(V)_{S^1}$ 
actually lies in $\ho^*_{S^1}(X)$.

The next result is used in the proof of 
Lemma~\ref{transfer-formula}.

\blemma \label{real-class}
Let $V$ be an orientable $S^1$-equivariant even dimensional
real vector bundle over $X$.  Suppose we are given two 
orientations of $V$,
which we denote by $V_{or_1}$ and $V_{or_2}$.  Define
$\sigma = 0$ if $V_{or_1} = V_{or_2}$, and $\sigma = 1$
otherwise.  Suppose $Q(x)$ is a power series such that
$Q(-x)=\gg Q(x)$, where $\gg=\pm 1$.  Then
  $$\mu_Q(V_{or_1}) = \gg^\sigma \mu_Q(V_{or_2}) \period$$
\elemma

\bp
\bi
\item[a)] If $Q(-x)=Q(x)$, $\mu_Q(V)$ is a power series in 
the equivariant Pontrjagin classes $p_j(V)_{S^1}$.  But 
Pontrjagin classes are independent of the orientation, so
$\mu_Q(V_{or_1}) = \mu_Q(V_{or_2})$.
\item[b)] If $Q(-x)=-Q(x)$, then $Q(x)=x\tilde{Q}(x)$, with
$\tilde{Q}(-x)=\tilde{Q}(x)$.  Hence 
$\mu_Q(V)=e_{S^1}(V)\cdot \mu_{\tilde{Q}}(V)$.
$e(V)_{S^1}$ changes sign when orientation changes sign, 
while $\mu_{\tilde{Q}}(V)$ is invariant, because of a).
\ei
\ep

\subsection{Acknowledgements}

I thank Matthew Ando for suggesting that I study the 
relationship between rigidity and Thom classes in equivariant 
elliptic cohomology.  I am also indebted to Mike Hopkins,
Jack Morava, and an anonymous referee for helpful comments.
Most of all I thank my advisor, Haynes Miller, who started
me on this subject, and gave me constant guidance and support.

\vspace{3mm}
\noindent
\textsc{department of mathematics, m.i.t., cambridge, ma 02139}\\
\textit{E-mail address:} \texttt{ioanid@math.mit.edu}

\end{document}